\documentclass{amsart}
\usepackage{amsmath}
\usepackage{amstext}
\usepackage{amssymb}
\usepackage{amsthm}
\swapnumbers
\theoremstyle{plain}
\newtheorem{thm}{Theorem}[section]
\newtheorem{lem}[thm]{Lemma}
\newtheorem{prop}[thm]{Proposition}
\newtheorem{cor}[thm]{Corollary}

\theoremstyle{definition}
\newtheorem{defi}[thm]{Definition}
\newtheorem{defs}[thm]{Definitions}
\newtheorem{defrmks}[thm]{Definition and Remarks}
\newtheorem{ex}[thm]{Example}

\theoremstyle{remark}
\newtheorem*{note}{Note}
\newtheorem{rmk}[thm]{Remark}

\DeclareMathOperator{\Ker}{Ker}

\DeclareMathOperator{\Ass}{Ass} \DeclareMathOperator{\ass}{ass}
\DeclareMathOperator{\red}{red}

\DeclareMathOperator{\Spec}{Spec}

\def\Z{\mathbb Z}
\def\N{\mathbb N}

\def\a{{\mathfrak{a}}}
\def\fa{{\mathfrak{a}}}
\def\fA{{\mathfrak{A}}}
\def\fb{{\mathfrak{b}}}

\def\fd{{\mathfrak{d}}}

\def\fm{{\mathfrak{m}}}
\def\fM{{\mathfrak{M}}}

\def\fp{{\mathfrak{p}}}
\def\fP{{\mathfrak{P}}}
\def\fq{{\mathfrak{q}}}

\def\fQ{{\mathfrak{Q}}}
\def\fr{{\mathfrak{r}}}
\def\fs{{\mathfrak{s}}}

\def\R{R_{\red}}
\def\nn{\relax\ifmmode{\mathbb N_{0}}\else$\mathbb N_{0}$\fi}
\def\lra{\longrightarrow}

\begin{document}

\title[IDEALS IN A PERFECT CLOSURE]{IDEALS IN A PERFECT CLOSURE,
LINEAR GROWTH OF PRIMARY DECOMPOSITIONS, AND TIGHT CLOSURE}
\author{RODNEY Y. SHARP}
\address[Sharp]{Department of Pure Mathematics,
University of Sheffield, Hicks Building, Sheffield S3 7RH, United Kingdom\\
{\it Fax number}: 0044-114-222-3769}
\email{R.Y.Sharp@sheffield.ac.uk}
\author{NICOLE NOSSEM}
\address[Nossem]{Department of Pure Mathematics,
University of Sheffield, Hicks Building, Sheffield S3 7RH, United Kingdom\\
{\it Fax number}: 0044-114-222-3769}
\email{N.Nossem@sheffield.ac.uk}
\thanks{The first author was partially supported by the
Engineering and Physical Sciences Research Council of the United
Kingdom (Overseas Travel Grant Number GR/S11459/01), and the
Mathematical Sciences Research Institute, Berkeley.
The second author was supported by a fees-only studentship provided
by the Engineering and Physical Sciences Research Council of the
United Kingdom.}


\date{\today}


\begin{abstract}
This paper is concerned with tight closure in a commutative
Noetherian ring $R$ of prime characteristic $p$, and is motivated
by an argument of K. E. Smith and I. Swanson that shows that, if
the sequence of Frobenius powers of a proper ideal $\fa$ of $R$
has linear growth of primary decompositions, then tight closure
(of $\fa$) `commutes with localization at the powers of a single
element'.  It is shown in this paper that, provided $R$ has a
weak test element, linear growth of primary
decompositions for other sequences of ideals of $R$ that
approximate, in a certain sense, the sequence of Frobenius powers
of $\fa$ would not only be just as good in this context, but,
in the presence of a certain additional finiteness property, would
actually imply that tight closure (of $\fa$) commutes with
localization at an arbitrary multiplicatively closed subset of
$R$.

Work of M. Katzman on the localization problem for tight closure raised the
question as to whether the union
of the associated primes of the tight closures of the Frobenius powers of $\fa$
has only finitely many maximal members. This paper develops, through
a careful analysis of the ideal
theory of the perfect closure of $R$, strategies
for showing that
tight closure (of a specified ideal $\fa$ of $R$)
commutes with localization at an arbitrary multiplicatively closed subset of $R$
and for showing that the union
of the associated primes of the tight closures of the Frobenius powers of $\fa$
is actually a finite set. Several applications of the strategies are
presented; in most of them it was already known that tight closure commutes
with localization, but the resulting affirmative answers to Katzman's question in the
various situations considered are believed to be new.
\end{abstract}

\maketitle

\setcounter{section}{-1}
\section{\sc Introduction}
\label{in}

This paper was motivated by a desire to explore the property,
which might be possessed by certain sequences of proper ideals in
a commutative Noetherian ring of prime characteristic, of having linear
growth of primary decompositions. The property has been studied by
K. E. Smith and I. Swanson \cite{KESISLinGrowth} in the context of the
localization problem for tight closure.

To set the scene, let $R$ be a commutative Noetherian ring of
prime characteristic $p$, and let $\fa$ be a proper ideal of $R$.
For $n \in \nn$ (we use $\nn$ (respectively $\N$) to denote the
set of non-negative (respectively positive) integers), the {\em
$n$-th Frobenius power\/} $\fa^{[p^n]}$ of $\fa$ is the ideal of
$R$ generated by all $p^n$-th powers of elements of $\fa$. Also
$R^{\circ}$ denotes the complement in $R$ of the union of the
minimal prime ideals of $R$.

An element $r \in R$ belongs to the {\em tight closure $\fa^*$ of
$\fa$\/} if and only if there exists $c \in R^{\circ}$ such that
$cr^{p^n} \in \fa^{[p^n]}$ for all $n \gg 0$. The theory of tight
closure was invented by M. Hochster and C. Huneke \cite{HocHun90},
and many applications have been found for the theory: see
\cite{Hunek96}. A major open problem about this theory concerns
the question of whether tight closure `commutes with
localization': see \cite[Chapter 12]{Hunek96}. It is known that it
does for many particular choices of $R$ and $\fa$: see
Aberbach--Hochster--Huneke \cite{AHH}, and also \cite[Chapter
12]{Hunek96}.

The starting point for the work in this paper is the following: if
the Frobenius powers of the proper ideal $\fa$ of $R$ have {\em
linear growth of primary decompositions\/} in the sense that there
exists a positive integer $h$ such that, for every non-negative
integer $n$, there exists a minimal primary decomposition
$\fa^{[p^n]} = \fq_{1,n} \cap \ldots \cap \fq_{k_n,n}$ with
$\sqrt{\fq_{i,n}}^{[p^n]h} \subseteq \fq_{i,n}$ for all $i =1,
\ldots, k_n$, then tight closure of $\fa$ commutes with
localization at a multiplicatively closed subset consisting of the
powers of a single element of $R$, that is, $\fa^*R_u = (\fa R_u)^*$
 for all $u \in R$.
A proof of this fact can be extracted
from the proof of \cite[Corollary (1.3)]{KESISLinGrowth}.

It should be noted that Swanson established an analogous linear growth property
for the ordinary powers of a proper ideal in an arbitrary
commutative Noetherian ring in \cite[Theorem 3.4]{ISPowerIdeal}.
The present first author subsequently provided a shorter proof of
a more general result in \cite{67}.

In this paper, we shall approach questions about linear growth of
primary decompositions by studying the ideals in the perfect
closure $R^{\infty}$ of $R$. When $R$ is not reduced, by the
perfect closure of $R$ we mean the perfect closure of $\R :=
R/\sqrt{(0)}$: we recall M. J. Greenberg's work \cite{Green65} in
this context in \S\ref{pc}.

A study of the ideals of $R^{\infty}$ leads us to the concept of
an {\em $f$-sequence of ideals\/} of $R$, that is, a sequence
$(\fa_n)_{n\in \nn}$ of ideals of $R$ such that $f^{-1}(\fa_{n+1})
= \fa_n$ for all $n \in\nn$, where $f : R \lra R$ denotes the
Frobenius homomorphism. The sequence $((\fa^{[p^n]})^*)_{n\in
\nn}$ of tight closures of the Frobenius powers of the ideal $\fa$
of $R$ is one example of an $f$-sequence. Another comes from
taking the $F$-closures of the Frobenius powers of $\fa$: the
$F$-closure $\fa^F$ of $\fa$ is defined as $$ \fa ^F := \left\{ r
\in R : \mbox{~there exists~} n \in \nn \mbox{~such that~} r^{p^n}
\in \fa^{[p^n]}\right\},
$$
and  $((\fa^{[p^n]})^F)_{n\in \nn}$ is another example of an
$f$-sequence. A third example involves the so-called plus closures
in the case when $R$ is a domain: the {\em plus closure\/} $\fa^+$
of $\fa$ is defined to be the contraction back to $R$ of the
extension of $\fa$ to the integral closure of $R$ in an algebraic
closure of its field of fractions, and we shall see in \S\ref{efs}
that $((\fa^{[p^n]})^+)_{n\in \nn}$ is an $f$-sequence.

It turns out that the set of $f$-sequences of ideals of $R$ is in
a natural bijective correspondence with the set of ideals of the
perfect closure $R^{\infty}$ of $R$. We explore this bijective
correspondence in some detail in the early part (\S\ref{jp},
\S\ref{pd} and \S\ref{pdpc}) of this paper, building on work of D.
A. Jordan \cite{Jorda82} that provides, among other things, a
rather concrete description of $R^{\infty}$.

The sequence of Frobenius powers of an ideal $\fa$ of $R$ is not
always an $f$-sequence. In \S\ref{lg} of the paper, we extend the
concept of linear growth of primary decompositions to
$f$-sequences in the obvious way. One of the main results of that
section is Theorem \ref{pd.7}, which shows that if a proper ideal
of $R^{\infty}$ has a primary decomposition, then the $f$-sequence
$(\fa_n)_{n\in \nn}$ of ideals of $R$ to which it corresponds has
linear growth of primary decompositions, and $\bigcup_{n \in \nn}
\ass \fa_n$ is finite. It should be noted that $R^{\infty}$ is
only Noetherian in rather uninteresting situations, so that the
existence of a primary decomposition for a proper ideal of
$R^{\infty}$ would be a bonus that we should not expect in all
cases.

Several of the main results of this paper involve the hypothesis
that $R$ has a $p^{m_0}$-weak test element for some $m_0 \in \nn$.
A {\em $p^{m_0}$-weak test element\/} is an element $c' \in
R^{\circ}$ such that, for every ideal $\fb$ of $R$ and for $r \in
R$, it is the case that $r \in \fb^*$ if and only if $c'r^{p^n}
\in \fb^{[p^n]}$ for all $n \geq m_0$. A $p^0$-weak test element
is called a {\em test element\/}. It is a result of Hochster and
Huneke \cite[Theorem (6.1)(b)]{HocHun94} that an algebra of finite
type over an excellent local ring of characteristic $p$ has a
$p^{m_0}$-weak test element for some $m_0 \in \nn$.

In Theorem \ref{lg.3}, we relate linear growth of primary
decompositions for certain $f$-sequences to the general
localization problem for tight closure, on the assumption that $R$
has a weak test element. This result, when used in conjunction
with the above-mentioned Theorem \ref{pd.7}, has the following
consequence: if $R$ has a $p^{m_0}$-weak test element, $\fa$ is a
proper ideal of $R$ and $(\fa_n)_{n\in \nn}$ is an $f$-sequence of
ideals of $R$ that approximates the Frobenius powers of $\fa$ in
the sense that $\fa^{[p^n]} \subseteq \fa_n \subseteq
(\fa^{[p^n]})^*$ for each $n \in \nn$, and if the ideal of
$R^{\infty}$ to which the $f$-sequence $(\fa_n)_{n\in \nn}$
corresponds has a primary decomposition, then $\fa^*S^{-1}R = (\fa
S^{-1}R)^*$ for every multiplicatively closed subset $S$ of $R$.
Furthermore, we show in Theorem \ref{lg.33} that, in these
circumstances, the ideal of $R^{\infty}$ to which the $f$-sequence
$\big((\fa^{[p^n]})^*\big)_{n \in \nn}$ corresponds also has a
primary decomposition, so that $\big((\fa^{[p^n]})^*\big)_{n \in
\nn}$ has linear growth of primary decompositions and the set
$\bigcup_{n \in \nn} \ass (\fa^{[p^n]})^*$ is finite. This is of
interest because M. Katzman's approach in \cite{Katzm96} to the
localization problem for tight closure led to the following
question: is it always the case that $\bigcup_{n \in \nn}\ass
(\fa^{[p^n]})^*$ has only finitely many maximal elements?

These results in \S\ref{lg} together provide a strategy (presented
in detail in Theorem \ref{lg.3a}) for attempting to show that
tight closure of a given proper ideal $\fa$ of $R$ commutes with
localization and attempting to answer Katzman's question for
$\fa$. A variant of the strategy that applies when the Frobenius
powers of $\fa$ have linear growth of primary decompositions and
$\bigcup_{n \in \nn} \ass \fa^{[p^n]}$ is finite is presented in
Theorem \ref{lg.4}, although this should be accompanied by the
warning that Katzman \cite{Katzm96} has provided an example of a
proper ideal $\fd$ in a ring of the type under consideration in
this paper for which $\bigcup_{n \in \nn} \ass \fd^{[p^n]}$ is
infinite!

Several applications of these strategies are presented in the
final \S\ref{ex}. We point out now that in most of those
applications, it was already known that tight closure commutes
with localization, often from the paper of
Aberbach--Hochster--Huneke \cite{AHH}; however, the numerous
results from use of the strategies that show that, in various
circumstances, the union of the associated primes of the tight
closures of the Frobenius powers of a fixed ideal is a finite set
are, we believe, new, and, in view of the above-mentioned question
of Katzman, of some interest.

For some of the applications, primary decompositions of appropriate ideals
of perfect closures are constructed by `approximating', in some sense,
the original ring $R$ by a regular commutative Noetherian ring of characteristic
$p$. Here, we just mention two of the applications, to give the general flavour.
One of the main results of \S\ref{ex} is
that, if $R$ is equidimensional and integral over a regular excellent subring
$A$, and $\fa$ is the extension to $R$ of a proper
ideal of $A$, then, provided $R$ has
a weak test element, the set $\bigcup_{n \in \nn} \ass (\fa^{[p^n]})^*$ is finite.
A corollary is that, if $x_1, \ldots, x_d$ form a system of parameters for an
equidimensional excellent local ring $R$, and
$\fa$ is a proper ideal of $R$ generated by `polynomials' in
$x_1, \ldots, x_d$ with coefficients in the prime subfield of $R$, then
the set $\bigcup_{n \in \nn} \ass (\fa^{[p^n]})^*$ is finite.

\section{\sc The perfect closure}
\label{pc}

For a commutative ring $R$ of prime characteristic $p$ we shall
always denote by $f:R\lra R$ the Frobenius homomorphism, for which
$f(r) = r^p$ for all $r \in R$. Recall that $R$ is said to be {\em
perfect\/} if $f$ is an isomorphism.

\begin{defi}
\label{pc.1}
M. J. Greenberg \cite[\S2]{Green65} proved that, for
a general, not necessarily perfect, $R$ of the above type,
there exists a pair
 $(R^\infty, \phi)$ where $R^\infty$ is a perfect ring and
 $\phi : R \lra R^\infty$ is a ring homomorphism such that, for
 any other such pair $(R', \phi')\/,$ there exists a unique
 ring homomorphism $\psi : R^\infty \lra R'$ with $\psi \circ \phi =
 \phi'$; furthermore, $\Ker \phi = \sqrt{(0)}$, the nilradical
 of $R$, and if $\alpha \in R^\infty\/,$ then $\alpha^{p^m} \in
 \phi(R)$ for some $m \geq 0\/.$ The ring $R^\infty$ is referred to
as the {\em perfect closure\/} of $R$.
\end{defi}

\begin{rmk}
\label{pc.2} It follows from Greenberg's work cited in \ref{pc.1}
that, with the notation of that definition, $R^\infty$ can be
constructed by forming the perfect closure of the reduced ring $\R
:= R/\sqrt{(0)}$. In \S\ref{jp} below, we shall provide a concrete
method of construction of the perfect closure of a reduced
commutative ring of prime characteristic.
\end{rmk}

\begin{rmk}
\label{pc.3} It also follows from Greenberg's work cited in \ref{pc.1} that,
if $R \cong \prod_{i=1}^n R_i$ is a finite product
 of commutative rings of the same prime characteristic,
 then $R^\infty \cong \prod_{i=1}^n {R_i}^\infty\/.$
\end{rmk}

\begin{rmk}
\label{ex.6} Let $A$ be a subring of a commutative ring $A'$ of
prime characteristic $p$. Then $A'$ is said to be a \emph{purely
inseparable extension of} $A$ if, for every element $a' \in A'\/,$
there exists
 $n = n(a') \in \nn$ such that
 ${a'}^{p^n} \in A$.  It is easy to deduce from \ref{pc.1} that
 this is the case if and only if the ring homomorphism
 $\iota^{\infty} : A^{\infty} \lra A'^{\infty}$ induced by the inclusion
 homomorphism $\iota : A \stackrel{\subseteq}{\lra} A'$ is an isomorphism.
 \end{rmk}

\section{\sc Jordan's construction}
\label{daj}

Throughout this section, let $A$ denote a (not necessarily commutative)
ring (with identity), and let $g : A \lra A$ be an injective ring
homomorphism. In \cite{Jorda82}, D. A. Jordan constructed in this
situation a ring $A'$ containing $A$
as a subring such that $g$ extends to $A'$ and is an automorphism of
$A'$; the same paper provides a detailed description of the left ideals of
$A'$. We now recall his construction of $A'$.

\begin{defi}
\label{daj.1} The
 {\em skew polynomial ring $A[x,g]$ associated to $A$ and $g$}
in an indeterminate
$x$ over $A$ is, as a left $A$-module, freely generated by $(x^i)_{i \in \nn}$,
 and so consists
 of all polynomials $\sum_{i = 0}^n a_i x^i$, where  $n \in \nn$
 and  $a_0,\ldots,a_n \in A$; however, its multiplication is subject to the
 rule
 $$
  xa = g(a)x \quad \mbox{~for all~} a \in A\/.
 $$
 The {\em (generalized) skew Laurent polynomial ring $A[x,x^{-1},g]$
associated to $A$ and
 $g$} is the left quotient ring of $A[x,g]$ with respect to the left
 denominator set
 $\{x^i : i \in \nn\}.$ Hence elements in $A[x,x^{-1},g]$ are finite
 sums of elements of the form $x^{-j}ax^i$ where $a \in A$ and
 $i, j \in \nn$ and multiplication in $A[x, x^{-1},g]$ is subject to the
 rules
 $$
  xa = g(a)x, \quad ax^{-1} = x^{-1}g(a) \quad \mbox{~for all~} a \in A.
 $$
(In the case when $g$ is an automorphism, the elements of $A[x,x^{-1},g]$
 can be written in the form $\sum_{i = m}^n a_ix^i$ where
 $m, n \in \Z$ and
 $a_m, \ldots, a_n \in A$. Then $A[x,x^{-1},g]$
 is the standard skew Laurent polynomial ring.)

 Let $A'$ be the subring of the skew Laurent
 polynomial ring $A[x,x^{-1},g]$ consisting of all
 elements of the form $x^{-i}ax^i$ with $i \in \nn$ and $a \in A$.
 As $x^{-i}ax^i = x^{-(i+j)}g^j(a)x^{i+j}$ for all $j \in \nn,$ this is a
 subring of $A[x,x^{-1},g]$ with
 $$
  x^{-i}ax^i + x^{-j}bx^j = x^{-(i+j)}(g^j(a)+g^i(b))x^{i+j}, \quad
(x^{-i}ax^i) ( x^{-j}bx^j) = x^{-(i+j)}(g^j(a)g^i(b))x^{i+j}
   $$
   for all $i,j \in \nn$ and $a,b \in A$.
 The ring homomorphism $g$ extends to $A'$ via
 $g(x^{-i}ax^i) = x^{-i}g(a)x^i$ for all $i \in \nn$
 and $a \in A$; in fact, $g : A' \lra A'$ is a ring isomorphism.
\end{defi}

Jordan's philosophy in \cite{Jorda82} is based on the fact that
the skew Laurent polynomial rings $A[x,x^{-1},g]$ and
$A'[x,x^{-1},g]$ are isomorphic; this means that one may tackle
problems concerning the structure of $A[x,x^{-1},g]$ by reducing
to the more familiar case where $g$ is an automorphism, provided
that one can handle the relationship between $A$ and $A'$. To this
end, Jordan provided in \cite[\S4]{Jorda82} a comprehensive
description of the left ideals of $A'$. We review this
next.

\begin{defs}
\label{daj.2} Let the situation be as in \ref{daj.1}.

\begin{enumerate}
\item Jordan calls a left ideal $I$ of $A$ {\em closed\/} if
 $g^{-n}(Ag^n(I))~\left(:= (g^{n})^{-1}(Ag^n(I))\right) \subseteq I$
  for all $n \in \N.$

  The argument in the proof of \cite[Proposition 4.4(i)]{Jorda82}
  shows that, if $J$ is an arbitrary left ideal of $A$ and $k \in
  \nn$, then $\rho_k(J) := \bigcup_{n \in \nn} g^{-n}(Ag^{n+k}(J))$
  is always a closed left ideal of $A$.

\item Jordan says that a sequence $(I_n)_{n \in \nn}$
of left ideals of $A$ is a {\em $g$-sequence\/} if
$g^{-1}(I_{n+1}) = I_n$ for all $n \in \nn$. (Actually, Jordan
also required that all the left ideals $I_n$ in the sequence be
closed; however, this is an automatic consequence of the property
that $g^{-1}(I_{n+1}) = I_n$ for all $n \in \nn$, as we now show.
Let $n \in \nn$ and $i \in \N$. Since $g^{-i}(I_{n+i}) = I_n$, we
have $Ag^i(I_n) \subseteq I_{n+i}$. Therefore $g^{-i}(Ag^i(I_n))
\subseteq g^{-i}(I_{n+i}) = I_n$. This shows that $I_n$ is
closed.)

\item If $J$ is a left ideal of $A$, then the argument in the proof
of \cite[Proposition 4.4(ii)]{Jorda82} shows that
$\left(\rho_n(J)\right)_{n \in \nn}$ is a $g$-sequence: we refer
to this as the {\em canonical $g$-sequence associated to $J$}.

\item Jordan defined \cite[p.\ 439]{Jorda82} a natural partial order
on the set of $g$-sequences of left ideals of $A$: for such
$g$-sequences $\left(I_n\right)_{n \in \nn}$ and
$\left(J_n\right)_{n \in \nn}$, we write
$$
\left(I_n\right)_{n \in \nn} \leq \left(J_n\right)_{n \in \nn}
\quad \mbox{if and only if} \quad I_n \subseteq J_n \mbox{~for
all~} n \in \nn.
$$

\item We shall occasionally find it convenient to use the ring
homomorphisms $\phi_n : A \lra A'~(n \in \nn)$ defined by
$\phi_n(a) = x^{-n}ax^n$ for all $a \in A$ (and $n \in \nn$).
\end{enumerate}
\end{defs}

\begin{thm}
\label{daj.3} {\rm (Jordan)} Let the situation be as in\/ {\rm
\ref{daj.1}}.

There is an order-preserving bijection, $\Gamma\/,$ from the
set of left
 ideals of $A'$, partially ordered by inclusion, to the partially-ordered set
 of $g$-sequences of left ideals of $A$ given by
 $$
  \Gamma : I \longmapsto (I_n)_{n \in \nn} \quad \mbox{~where~}
     I_n := \{ a \in A : x^{-n}ax^n \in I\} = \phi_n^{-1}(I)
     \mbox{~for all~} n \in \nn.
 $$
 The inverse bijection, $\Gamma^{-1},$ also order-preserving, is given by
 $$
\Gamma^{-1} : (I_n)_{n \in \nn} \longmapsto \bigcup_{n \in \nn} x^{-n}I_nx^n =
  \bigcup_{n \in \nn} \phi_n(I_n).
 $$
\end{thm}

If, with the notation of the above theorem,
$\Gamma(I) = (I_n)_{n \in \nn}\/,$ then we shall say that
$(I_n)_{n \in \nn}$ is the
{\em $g$-sequence associated to, or corresponding to, $I$\/}, and,
given a $g$-sequence $(I_n)_{n \in \nn},$ we shall call the left ideal
$\Gamma^{-1}((I_n)_{n \in \nn})$ the {\em left ideal associated to, or
corresponding to, the $g$-sequence $(I_n)_{n \in \nn}.$}

\section{\sc Application of Jordan's construction to perfect closures}
\label{jp}

In this section, we explain the relevance of Jordan's construction
of \S \ref{daj} to the perfect closure of a reduced commutative
Noetherian ring $R$ of prime characteristic $p$ (and $R$ will have
this meaning throughout this section). Since $R$ is reduced, the
Frobenius homomorphism $f : R \lra R$ is injective, and we may
take $f : R \lra R$ for the $g : A \lra A$ in Jordan's
construction of \S \ref{daj}.

\begin{rmk}
\label{jp.1} We apply Jordan's construction with $f : R \lra R$ in
the r\^ole of $g : A \lra A$.
\begin{enumerate}
\item The ring $A'$ (which is commutative in this case)
is just the perfect closure of $R$, because if
$R'$ is a perfect commutative ring and
 $\phi' : R \lra R'$ is a ring homomorphism, then there is a unique
 ring homomorphism $\psi : A' \lra R'$ whose restriction to $R$ is
 $\phi'$, given by $\psi(x^{-n}rx^n) = \phi'(r)^{1/p^n}$,
 the unique $p^n$-th root of $\phi'(r)$ in $R'$, for all $n \in \nn$ and
 $r \in R$.

We shall denote $A'$ in this case by $R^{\infty}$. Observe that,
for $n \in \nn$ and $r \in R$, the element $x^{-n}rx^n$ of
$R^{\infty}$ is the unique $p^n$-th root $r^{1/p^n}$ of $r$ in
$R^{\infty}$, and the ring homomorphism $\phi_n : R \lra
R^{\infty}$ of \ref{daj.2}(v) maps $r$ to $r^{1/p^n}$. Thus
Jordan's construction provides a rather concrete presentation of
the perfect closure of $R$.

\item A sequence $(\a_n)_{n \in \nn}$
of ideals of $R$ is an {\em $f$-sequence\/} if $f^{-1}(\a_{n+1}) =
\a_n$ for all $n \in \nn$.

\end{enumerate}
\end{rmk}

Of course, Jordan's Theorem \ref{daj.3} provides detailed
information about the ideals of $R^{\infty}$.

\begin{cor}
\label{jp.2}
There is an order-preserving bijection,
$\Gamma\/,$ from the
set of
 ideals of $R^{\infty}$, partially ordered by inclusion,
 to the partially-ordered set
 of $f$-sequences of ideals of $R$ given by
 $$
  \Gamma : \fA \longmapsto (\a_n)_{n \in \nn} \quad \mbox{~where~}
     \a_n := \{ r \in R : x^{-n}rx^n = r^{1/p^n} \in \fA\} = \phi_n^{-1}(\fA)
     \mbox{~for all~} n \in \nn.
 $$
 The inverse bijection, $\Gamma^{-1},$ also order-preserving, is given by
 $$
  \Gamma^{-1} : (\a_n)_{n \in \nn} \longmapsto \bigcup_{n \in \nn}
  x^{-n}\a_nx^n =   \bigcup_{n \in \nn} \phi_n(\a_n).
 $$
\end{cor}

\section{\sc Properties of $f$-sequences}
\label{pd}

Throughout this section and the remainder of the paper, $R$ will
denote a commutative Noetherian ring of prime characteristic
$p$\/; note that we have dropped the hypothesis that $R$ be
reduced. In this section, we extend the concept of $f$-sequence of
ideals to this more general situation, and then develop the
concept in some detail.

\begin{defs}
\label{pd.0} Let $\fa$ and $\fb$ be ideals of $R$.

\begin{enumerate}
\item The ideal $\a$ of $R$ is said to be {\em $F$-closed\/} if, whenever
$r \in R$ is such that $r^{p^n} \in \a^{[p^n]}$ for some $n \in
\nn$ (or, equivalently, for all $n \gg 0$), then $r \in \a$.

\item An {\em $f$-sequence of ideals of $R$\/} is a sequence
$(\a_n)_{n \in \nn}$ of ideals of $R$ such that $f^{-1}(\a_{n+1})
= \a_n$ for all $n \in \nn$. There is a partial order on the set
of such $f$-sequences, defined analogously to the partial order of
\ref{daj.2}(iv). The argument of \ref{daj.2}(ii) applies here to
show that every term in an $f$-sequence is $F$-closed.

\item The {\em $F$-closure of $\fb$\/} is the
($F$-closed) ideal $\fb^F := \bigcup_{n \in \nn}
f^{-n}(Rf^{n}(\fb))= \bigcup_{n \in \nn}
f^{-n}\big(\fb^{[p^n]}\big)$, and $\big((\fb^{[p^n]})^F\big)_{n
\in \nn}$ is an $f$-sequence, called the {\em canonical
$f$-sequence associated to $\fb$\/}.
\end{enumerate}
\end{defs}

\begin{rmk}
\label{pd.1a} Let $(\a_n)_{n \in \nn}$ be an $f$-sequence of
ideals of $R$.

\begin{enumerate}
\item The $f$-sequence
 $(\a_n)_{n \in \nn}$ is descending, because
 if
 $a \in \a_{n + 1}$ for some $n \in \nn$, then
 $f(a) = a^p \in \a_{n+1}\/,$ so that $a \in f^{-1}(\a_{n+1}) =
 \a_n\/.$ Further, $\a_n^{[p]} \subseteq \a_{n+1}$ since $f(\a_n)
 \subseteq \a_{n+1}\/.$

 This implies, as $\a_{n+1}$ is $F$-closed,
  that $\a_n^{[p]} \subseteq
 (\a_n^{[p]})^F \subseteq \a_{n+1} \subseteq \a_n$, so that
 $\sqrt{\a_n} = \sqrt{\a_{n+1}}$,
 for all $n \in \nn\/.$

Similar arguments show that $(\fa_n^{[p^j]})^F \subseteq \fa_{n+j}$ for all
$n,j \in \nn$.
\item Each term $\fa_n$ in the $f$-sequence contains $\sqrt{(0)}$,
since $\fa_n$ is $F$-closed.
\item We shall sometimes need to use the fact that if, for some $t \in
\N$, one has a sequence $(\fb_n)_{n \geq t}$ of ideals of $R$ such
that $f^{-1}(\fb_{n+1}) = \fb_n$ for all $n \geq t$, then there is
exactly one way of extending the sequence `downwards' to an
$f$-sequence $(\fb_n)_{n \in \nn}$, and one achieves this
extension by setting $\fb_n = f^{-(t-n)}(\fb_t)$ for all $n = t-1,
\ldots,0$.
\item Note that $\ass \a_n \subseteq \ass \a_{n+1}$ for each $n \in \nn$, since
$f^{-1}(\a_{n+1}) = \a_n$.
\item Note also that all the terms of the $f$-sequence $(\a_n)_{n \in
\nn}$ have the same set of minimal primes, since $\sqrt{\a_n} =
\sqrt{\a_{n+1}}$ for all $n \in \nn$, by (i). The members of the
common set of minimal primes of the terms of the $f$-sequence are
referred to as the \emph{minimal primes of the $f$-sequence}; note
that there are only finitely many of these minimal primes.
\end{enumerate}
\end{rmk}

The fact that every term in an $f$-sequence of ideals of $R$
contains the nilradical of $R$ (by \ref{pd.1a}(ii)) means that
there is an obvious bijective correspondence (given by extension
of the terms and contraction of the terms) between the set of
$f$-sequences of ideals of $R$ and the set of $f$-sequences of
ideals of $\R$. This means that Jordan's Corollary \ref{jp.2} has
an analogue that applies in this more general situation.

\begin{cor}
\label{pd.1j} We shall use overlines to denote natural images in
$\R$ of elements of $R$ under the natural homomorphism $\pi : R
\lra \R$, although we shall sometimes omit the overlines when $R$
is itself reduced; we shall interpret the subring $A'$ of the skew
Laurent polynomial ring $\R[x,x^{-1},f]$, as defined in\/ {\rm
\ref{daj.1}}, as the perfect closure $(\R)^{\infty} = R^{\infty}$
of $R$; for each $n \in \nn$, we shall use $\phi_n : \R \lra
R^{\infty}$ to denote the ring homomorphism for which
$\phi_n(\overline{r}) = x^{-n}\overline{r}x^n
=\overline{r}^{1/p^n}$ for all $r \in R$, and $\psi_n :R \lra
R^{\infty}$ to denote the composition $\phi_n \circ \pi$.

There is an order-preserving bijection, $\Gamma\/,$ from the
set of
 ideals of $R^{\infty}$, partially ordered by inclusion,
 to the partially-ordered set
 of $f$-sequences of ideals of $R$ given by
 $$
  \Gamma : \fA \longmapsto (\a_n)_{n \in \nn} \quad \mbox{~where~}
     \a_n := \{ r \in R : x^{-n}\overline{r}x^n =
     \overline{r}^{1/p^n} \in \fA\} = \psi_n^{-1}(\fA)
     \mbox{~for all~} n \in \nn.
 $$
 The inverse bijection, $\Gamma^{-1},$ also order-preserving, is given by
 $$
  \Gamma^{-1} : (\a_n)_{n \in \nn} \longmapsto \bigcup_{n \in \nn}
  x^{-n}\pi(\a_n)x^n =  \bigcup_{n \in \nn}
   \psi_n(\a_n).
 $$
\end{cor}

\begin{lem}
\label{bk.4} Let $\eta : R \lra R'$ be a homomorphism of commutative
Noetherian rings of characteristic $p$, and let
$(\fa'_n)_{n \in\nn}$ be an $f$-sequence of ideals of $R'$ corresponding
to the ideal $\fA'$ of $R'^{\infty}$.
Let $\eta^{\infty} : R^{\infty} \lra R'^{\infty}$ be
the ring homomorphism induced by $\eta$ (see\/ {\rm \S1}).

Then the $f$-sequence of ideals of $R$ corresponding
to the ideal $(\eta^{\infty})^{-1}(\fA')$ of $R^{\infty}$ is
$\left(\eta^{-1}(\fa'_n)\right)_{n \in\nn}$.
\end{lem}

\begin{proof} Let $\psi_n :R \lra R^{\infty}$ be the ring homomorphism of
\ref{pd.1j}, and let $\psi'_n :R' \lra R'^{\infty}$ be the
corresponding ring homomorphism for $R'$. The claim follows easily
from the fact that $\eta^{\infty} \circ \psi_n = \psi'_n \circ
\eta$ for all $n \in \nn$.
\end{proof}

\section{\sc Examples of $f$-sequences}
\label{efs}

Throughout this section, $R$ will denote a commutative Noetherian
ring of prime characteristic $p$. The first two lemmas in this
section present important examples of $f$-sequences.

\begin{lem}
\label{pd.1c} Let $\fa$ be an ideal of $R$. Then
$\big((\fa^{[p^n]})^*\big)_{n \in \nn}$, the sequence of tight
closures of the Frobenius powers of $\fa$, is an $f$-sequence.
\end{lem}

\begin{proof} Let $n \in \nn$.

Let $r \in f^{-1}\big((\fa^{[p^{n+1}]})^*\big)$, so that $r^p \in
(\fa^{[p^{n+1}]})^*$ and there exists $c \in R^{\circ}$ such that
$c(r^p)^{p^j} \in (\fa^{[p^{n+1}]})^{[p^j]}$ for all $j \gg 0$;
hence $cr^{p^{j+1}} \in (\fa^{[p^{n}]})^{[p^{j+1}]}$ for all $j \gg 0$,
so that $r \in (\fa^{[p^{n}]})^*$. Hence
$f^{-1}\big((\fa^{[p^{n+1}]})^*\big) \subseteq (\fa^{[p^{n}]})^*$.

Now let $r \in (\fa^{[p^{n}]})^*$, so that there exists $c \in
R^{\circ}$ such that $cr^{p^j} \in (\fa^{[p^{n}]})^{[p^{j}]}$ for
all $j \gg 0$. Then $c^pr^{p^{j+1}} \in
((\fa^{[p^{n}]})^{[p^{j}]})^{[p]}$ for all $j \gg 0$. Hence
$c^p(r^p)^{p^j} \in (\fa^{[p^{n+1}]})^{[p^{j}]}$ for all $j \gg
0,$ and so $r^p \in (\fa^{[p^{n+1}]})^*$ because $c^p \in
R^{\circ}$.
\end{proof}

The next lemma concerns the case where $R$ is a domain and
involves the plus closure of an ideal $\fa$ of $R$: the reader was
reminded about this concept in the Introduction. In this
situation, we denote by $R^+$ the integral closure of $R$ in an
algebraic closure of its field of fractions, and refer to $R^+$ as
the {\em absolute integral closure of \/} $R$. Recall from Huneke
\cite[p.\ 15]{Hunek96} that $\fa \subseteq \fa^+ \subseteq \fa^*$.

\begin{lem}
\label{fs.1} Assume that $R$ is a domain and let $\fa$ be an ideal
of $R$. Then $\big((\fa^{[p^n]})^+\big)_{n \in \nn}$, the sequence
of plus closures of the Frobenius powers of $\fa$, is an
$f$-sequence. In fact, it is the $f$-sequence corresponding to the
ideal $\fa R^+\cap R^{\infty}$ of $R^{\infty}$ in the
correspondence of\/ {\rm \ref{jp.2}}.
\end{lem}

\begin{proof} It is enough to prove the second statement. Let
$(\a_n)_{n \in \nn}$ be
the $f$-sequence of ideals of $R$ that corresponds to the ideal
$\fa R^+\cap R^{\infty}$ of $R^{\infty}$. Let $n \in \nn$, and
suppose that $a_1, \ldots, a_t$ generate $\fa$. Note that $\fa_n
=\{ r \in R : r^{1/p^n} \in \fa R^+\cap R^{\infty}\}. $

Let $r \in (\fa ^{[p^{n}]})^+$. Thus $r \in R$ and there exist
$\sigma_1, \ldots, \sigma_t \in R^+$ such that
$
r = a_1^{p^n}\sigma_1 + \cdots + a_t^{p^n}\sigma_t.
$
 Now, for each $i
= 1, \ldots, t$, the unique $p^n$-th root $\sigma_i^{1/p^n}$ of
$\sigma_i$ in the algebraic closure of the field of fractions of
$R$ is integral over $R^+$, and so belongs to $R^+$. Hence
$r^{1/p^n} = a_1\sigma_1 ^{1/p^n}+ \cdots + a_t\sigma_t ^{1/p^n}
\in \fa R^+ \cap R^{\infty}.
$
Therefore $r \in \fa_n$. Thus $(\fa ^{[p^{n}]})^+ \subseteq
\fa_n$. The reverse inclusion is even easier.
\end{proof}

We now drop the hypothesis that $R$ is a domain, and revert to our
standard hypotheses about $R$.
We would like to have a variant of the concept of plus closure available for use
in situations where $R$ is not a domain (and not even reduced), and the
following definition introduces a suitable one.

\begin{defrmks}
\label{me.13} Suppose that there are $h$ minimal prime ideals $\fp_1, \ldots, \fp_h$
of $R$; let $\fa$ be an ideal of $R$. For each $i = 1, \ldots, h$, one can use
the integral domain $R/\fp_i$ to construct ($(R/\fp_i)^+$ and) the plus closure
$((\fa + \fp_i)/\fp_i)^+$ of the ideal $(\fa + \fp_i)/\fp_i$ of $R/\fp_i$.

We define the {\em plus closure\/}
$\fa^+$ of $\fa$ to be the contraction back to $R$
of the ideal $\prod_{i=1}^h ((\fa + \fp_i)/\fp_i)^+$ of the direct product
$\prod_{i=1}^h R/\fp_i$ under the natural ring homomorphism
$\nu : R \lra \prod_{i=1}^h R/\fp_i$.

Note that, in view of \S1, one can identify
$\big(\prod_{i=1}^h R/\fp_i\big)^{\infty}$ with
$\prod_{i=1}^h (R/\fp_i)^{\infty}$ (under an isomorphism that maps
$\left((r_1 + \fp_1), \ldots, (r_h + \fp_h)\right)^{1/p^n}$
to $\left((r_1 + \fp_1)^{1/p^n}, \ldots, (r_h + \fp_h)^{1/p^n}\right)$
for $r_1, \ldots, r_h
\in R$ and $n \in \nn$), and there is an induced injective ring homomorphism
$$
\nu^{\infty}: R^{\infty} = (\R)^{\infty} \lra \prod_{i=1}^h (R/\fp_i)^{\infty} =
\bigg(\prod_{i=1}^h R/\fp_i\bigg)^{\infty}
$$
such that $\nu^{\infty} \left( \overline{r}^{1/p^n}\right) =
\left((r + \fp_1)^{1/p^n}, \ldots, (r + \fp_h)^{1/p^n}\right)$ for all $r \in R$
and $n \in \nn$. We use $\nu^{\infty}$ to identify $R^{\infty}$ as a subring of
$\prod_{i=1}^h (R/\fp_i)^{\infty}$; the latter ring is identified with a subring of
$\prod_{i=1}^h (R/\fp_i)^+$ in the obvious natural way. Note that the diagram
\[
\begin{picture}(300,75)(-150,-25)
\put(-110,40){\makebox(0,0){$ R $}} \put(-55,46){\makebox(0,0){$^{
\psi_0 }$}} \put(0,40){\makebox(0,0){$ R^{\infty} $}}
\put(-90,40){\vector(1,0){70}} \put(-106,10){\makebox(0,0)[l]{$^{
\nu }$}} \put(4,10){\makebox(0,0)[l]{$^{ \nu ^{\infty} }$}}
\put(-110,-20){\makebox(0,0){$ \prod_{i=1}^h (R/\fp_i) $}}
\put(-55,-12){\makebox(0,0){$^{ \prod_{i=1}^h \psi_0^{(i)} }$}}
\put(0,-20){\makebox(0,0){$ \prod_{i=1}^h (R/\fp_i)^{\infty} $}}
\put(55,-14){\makebox(0,0){$^{ \subseteq }$}}
\put(110,-20){\makebox(0,0){$ \prod_{i=1}^h (R/\fp_i)^+ $}}
\put(-75,-20){\vector(1,0){40}} \put(35,-20){\vector(1,0){40}}
\put(-110,30){\vector(0,-1){40}} \put(0,30){\vector(0,-1){40}}
\end{picture}
\]
(in which $\psi_0: R \lra R^{\infty}$ and $\psi_0^{(i)}: R/\fp_i
\lra (R/\fp_i)^{\infty}~ (i = 1, \ldots,h)$ are the natural ring
homomorphisms) commutes, and
\begin{align*}
\fa^+ &= \nu^{-1} \bigg(\prod_{i=1}^h ((\fa +
\fp_i)/\fp_i)^+\bigg) = \nu^{-1} \bigg( \Big(\prod_{i=1}^h
\psi_0^{(i)}\Big)^{-1} \Big( \prod_{i=1}^h ((\fa +
\fp_i)/\fp_i)(R/\fp_i)^+ \cap (R/\fp_i)^{\infty}\Big)
\bigg)\\
&= (\psi_0)^{-1} \bigg((\nu^{\infty})^{-1} \Big( \Big(
\prod_{i=1}^h ((\fa + \fp_i)/\fp_i)(R/\fp_i)^+\Big) \bigcap \Big(
\prod_{i=1}^h
(R/\fp_i)^{\infty}\Big) \Big) \bigg)\\
&= (\psi_0)^{-1} \bigg(\Big( \prod_{i=1}^h ((\fa +
\fp_i)/\fp_i)(R/\fp_i)^+\Big) \bigcap R^{\infty} \bigg).
\end{align*}
\end{defrmks}

When $R$ is a domain, $\fa^+$, as defined in \ref{me.13}, coincides
with the plus closure of $\fa$, as defined by Hochster and Huneke.
The next lemma demonstrates that, in the general
case, plus closure retains one important property of the Hochster--Huneke concept.

\begin{lem}
\label{me.14} Let $\fa$ be an ideal of $R$, and use the notation of\/
{\rm \ref{me.13}}.

\begin{enumerate}
\item We have $\fa \subseteq \fa^+ \subseteq \fa^*$.
\item The $f$-sequence of ideals of\/ $\prod_{i=1}^h R/\fp_i$ to which the
ideal\/ $\prod_{i=1}^h ((\fa + \fp_i)/\fp_i)(R/\fp_i)^+ \cap
(R/\fp_i)^{\infty}$ of\/ $\prod_{i=1}^h (R/\fp_i)^{\infty} =\big(
\prod_{i=1}^h R/\fp_i\big)^{\infty}$ corresponds is
$\left(\prod_{i=1}^h ((\fa^{[p^n]} + \fp_i)/\fp_i)^+\right)_{n \in
\nn}$.
\item Also, $\big((\fa^{[p^n]})^+\big)_{n \in \nn}$ is the $f$-sequence of ideals of
$R$ corresponding to the ideal
$$\Big(
\prod_{i=1}^h ((\fa + \fp_i)/\fp_i)(R/\fp_i)^+\Big) \bigcap R^{\infty}$$ of
$R^{\infty}$.
\end{enumerate}
\end{lem}

\begin{proof} (i) It is clear that
$\fa \subseteq \fa^+$. Let $r \in \fa^+$. Then
$$
\nu (r) = (r + \fp_1, \ldots, r + \fp_h) \in
\prod_{i=1}^h ((\fa + \fp_i)/\fp_i)^+ \subseteq
\prod_{i=1}^h ((\fa + \fp_i)/\fp_i)^* \subseteq \prod_{i=1}^h R/\fp_i.
$$
It then follows from \cite[Theorem 1.3(c)]{Hunek96} that $r \in \fa^*$.

(ii) This is an easy consequence of \ref{fs.1}.

(iii) This now follows easily from Lemma \ref{bk.4} applied to the ring
homomorphism $\nu : R \lra \prod_{i=1}^h R/\fp_i$.
\end{proof}

The next lemma will help us to construct new examples of
$f$-sequences from known ones.

\begin{lem}
\label{pd.2a} Let $(\fa_n)_{n \in\nn}$ be an $f$-sequence of
ideals of $R$ and let $\fA$ be the corresponding ideal of
$R^{\infty}$. Let $S$ be a multiplicatively closed subset
 of $R$, and, with the notation of\/ {\rm \ref{pd.1j}},
 set $T := \bigcup_{n \in \nn}x^{-n}\pi(S)x^n$, a multiplicatively
closed subset of $R^{\infty}$.
 Further let $\phantom{\fa}^e$ and $\phantom{\fa}^c$ stand for extension and
 contraction with respect to the natural ring homomorphism
$R \lra S^{-1}R$, and let $\phantom{\fa}^E$ and $\phantom{\fa}^C$
stand for extension and
 contraction with respect to the natural ring homomorphism
$R^{\infty} \lra T^{-1}(R^{\infty})$.

 Then $(\fa_n^{ec})_{n \in\nn}$ is also an $f$-sequence of ideals of $R$, and
the ideal of $R^{\infty}$ to which it corresponds is $\fA^{EC}$.
\end{lem}

\begin{proof}
It is straightforward to check that $T$ is a multiplicatively
closed subset of $R^{\infty}$. Thus $\fA^{EC}$ is an ideal of
$R^{\infty}$. Let $(\fa'_n)_{n \in\nn}$ be the corresponding
$f$-sequence of ideals of $R$. It is enough for us to show that
$\fa'_n = \fa_n^{ec}$ for each $n \in \nn$. Note that $\fa'_n =
\{a \in R : x^{-n}\overline{a}x^n \in \fA^{EC}\}$.

Let $a \in \fa_n^{ec}$. Thus there exists $s \in S$ such that $sa
\in \fa_n$. Hence $x^{-n}\overline{sa}x^n \in \fA$, so that
$x^{-n}\overline{s}x^nx^{-n}\overline{a}x^n \in \fA$. Since
$x^{-n}\overline{s}x^n \in T$, this means that
$x^{-n}\overline{a}x^n \in \fA^{EC}$ and $a \in \fa'_n$.

Now let $a \in \fa'_n$, so that $x^{-n}\overline{a}x^n \in
\fA^{EC}$ and there exists $\tau \in T$ such that $\tau
x^{-n}\overline{a}x^n \in \fA$. We can write $\tau =
x^{-k}\overline{s}x^k$ for some $k \in \nn$ and $s \in S$. If $k
\leq n$, then
$$
x^{-n}\overline{s^{p^{n-k}}a}x^n =
x^{-n}\overline{s}^{p^{n-k}}x^nx^{-n}\overline{a}x^n =
x^{-k}\overline{s}x^kx^{-n}\overline{a}x^n = \tau
x^{-n}\overline{a}x^n \in \fA.
$$
Hence, in this case, $s^{p^{n-k}}a \in \fa_n$, so that $a \in
\fa_n^{ec}$ because $s^{p^{n-k}} \in S$.

We now consider the case where $k > n$. Then
$$
x^{-n}\overline{sa}x^n =
x^{-n}\overline{s}x^nx^{-n}\overline{a}x^n =
x^{-k}\overline{s}^{p^{k-n}}x^kx^{-n}\overline{a}x^n =
x^{-k}\overline{s}^{p^{k-n}-1}x^k(\tau x^{-n}\overline{a}x^n) \in
\fA.
$$
Therefore $sa \in \fa_n$, so that $a \in \fa_n^{ec}$.
\end{proof}

\section{\sc Primary decompositions in perfect closures}
\label{pdpc}

Throughout this section, $R$ will denote a commutative Noetherian ring
of prime characteristic $p$. Our aim in the section is to establish necessary
and sufficient conditions for a proper ideal $\fA$ of the perfect closure
$R^{\infty}$ to have a primary decomposition; our conditions will
be phrased in terms of the $f$-sequence of ideals of $R$ that corresponds to
$\fA$.

\begin{thm}
\label{pd.2} Let $\fA$ be an ideal of $R^{\infty}$, and let
$(\a_n)_{n \in \nn}$ be the corresponding $f$-sequence of ideals
of $R$. We use the notation of \/ {\rm \ref{pd.1j}}.

\begin{enumerate}

\item The ideal $\fA$ is prime if and only if there is a $\fp \in \Spec (R)$
such that $\fp = \a_n$ for all $n \in \nn$. Thus the $f$-sequence
corresponding to a prime ideal $\fP$ of $R^{\infty}$ is the
constant $f$-sequence $\left(\pi^{-1}(\fP\cap \R)\right)_{n \in
\nn}$, and the correspondence of\/ {\rm \ref{pd.1j}} yields an
inclusion-preserving bijective correspondence between $\Spec
(R^{\infty})$ and $\Spec (R)$.

\item The sequence $\left(\sqrt{\a_n}\right)_{n \in \nn}$ is an
$f$-sequence, and it corresponds to the ideal $\sqrt{\fA}$ of
$R^{\infty}$. Thus $\sqrt{\bigcup_{n \in \nn} x^{-n}\pi(\a_n)x^n}
= \bigcup_{n \in \nn} x^{-n}\pi(\sqrt{\a_n})x^n$.

\item The ideal $\fA$ is radical if and only if there is a radical ideal
$\fb$ of $R$ such that $\fb = \a_n$ for all $n \in \nn$.

\item Let $\fP \in \Spec (R^{\infty})$, and let $\fp$ be the
corresponding prime ideal of $R$: see\/ {\rm (i)} above. Then
$\fA$ is a $\fP$-primary ideal of $R^{\infty}$ if and only if
$\a_n$ is $\fp$-primary for all $n \in \nn$.

\item Let $\left(\fA_{\lambda}\right)_{\lambda \in \Lambda}$ be a family
of ideals of $R^{\infty}$, and, for each $\lambda \in \Lambda$,
let $(\a_{\lambda,n})_{n \in \nn}$ be the $f$-sequence
corresponding to $\fA_{\lambda}$. Then $\left(\bigcap_{\lambda \in
\Lambda} \a_{\lambda,n}\right)_{n \in \nn}$ is the $f$-sequence of
ideals of $R$ corresponding to the ideal $\bigcap_{\lambda \in
\Lambda}\fA_{\lambda}$ of $R^{\infty}$. Consequently,
$$
\bigcup_{n \in \nn} \bigg(x^{-n}\pi\bigg(\bigcap_{\lambda \in
\Lambda} \a_{\lambda,n}\bigg)x^n\bigg) = \bigcap_{\lambda \in
\Lambda}\bigg(\bigcup_{n \in \nn} x^{-n}
\pi(\a_{\lambda,n})x^n\bigg).
$$
\end{enumerate}
\end{thm}

\begin{proof} Since each term in an $f$-sequence of ideals of $R$,
and each prime ideal of $R$, contains $\sqrt{(0)}$, it is enough for
us to prove this theorem under the additional assumption that $R$
is reduced; with this assumption, $\psi_n = \phi_n$ for all $n \in
\nn$.

(i) We have $f^{-1}(\a_{n+1}) = \a_n = \phi_n^{-1}(\fA)$ for all
$n \in \nn$, and it is immediate from this that if $\fA$ is prime,
then $(\a_n)_{n \in \nn}$ is constant with all its terms equal to
the same prime ideal.

Conversely, if there exists $\fp \in \Spec (R)$ such that $\a_n =
\fp$ for all $n \in \nn$, then it is routine to check that $\fA =
\bigcup_{n \in \nn} x^{-n}\fp x^n \in \Spec (R^{\infty})$: note
that $1$ cannot be written as $x^{-n}r x^n$ for any $r \in \fp$,
and that if $a,b \in R$ are such that $x^{-i}a x^ix^{-j}b x^j \in
\fA$ for some integers $i, j \in \nn$, then
$x^{-i-j}f^{j}(a)f^{i}(b) x^{i+j} \in \fA$, so that
$f^{j}(a)f^{i}(b) \in \fa_{i+j} = \fp$, from which it follows that
$a$ or $b$ lies in $\fp$.

(ii) Of course, $\sqrt{\fA}$ is an ideal of $R^{\infty}$; by
\ref{pd.1j}, the $f$-sequence to which it corresponds is
$$
\left(\phi_n^{-1}\left(\sqrt{\fA}\right)\right)_{n \in \nn} =
\left(\sqrt{\phi_n^{-1}(\fA)}\right)_{n \in \nn} =
\left(\sqrt{\a_n}\right)_{n \in \nn}.
$$

(iii) The ideal $\fA$ is radical if and only if $\fA=\sqrt{\fA}$;
by \ref{pd.1j}, this is the case if and only if the $f$-sequences
to which $\fA$ and $\sqrt{\fA}$ correspond are the same; and it
follows from part (ii) that this is the case if and only if $\a_n
= \sqrt{\a_n}$ for all $n \in \nn$. Finally, note that, when this
is the case, we have $\a_n = \a_{n+1}$ for each $n \in \nn$, by
\ref{pd.1a}(i).

(iv) If $\fA$ is $\fP$-primary, then, since $\fp = \phi_n^{-1}(\fP)$ (by
part (i)) and $\a_n = \phi_n^{-1}(\fA)$, it follows that $\a_n$ is
$\fp$-primary (for all $n \in \nn$).

Conversely, suppose that $\a_n$ is $\fp$-primary for all $i \in
\nn$. Then, by parts (i) and (ii), we have $\fA \subseteq
\sqrt{\fA} = \fP$, so that $\fA$ is proper. Suppose that $a,b \in
R$ are such that $x^{-i}a x^ix^{-j}b x^j \in \fA$ for some
integers $i, j \in \nn$. Then $x^{-i-j}f^{j}(a)f^{i}(b) x^{i+j}
\in \fA$, so that $f^{j}(a)f^{i}(b) \in \fa_{i+j}$. It follows
from this that either $f^{j}(a) \in \a_{i+j}$ or $f^{i}(b) \in
\sqrt{\a_{i+j}} = \fp$, that is, either $a \in f^{-j} (\a_{i+j}) =
\a_i$ or $b \in \fp$. Hence either $x^{-i}a x^i \in \fA$ or
$x^{-j}b x^j \in \fP = \sqrt{\fA}$. Thus $\fA$ is $\fP$-primary.

(v) The $f$-sequence corresponding to the ideal $\bigcap_{\lambda
\in \Lambda}\fA_{\lambda}$ of $R^{\infty}$ is, by \ref{pd.1j},
$$
\bigg(\phi_n^{-1}\bigg(\bigcap_{\lambda \in
\Lambda}\fA_{\lambda}\bigg)\bigg)_{n \in \nn} =
\bigg(\bigcap_{\lambda \in
\Lambda}\phi_n^{-1}\left(\fA_{\lambda}\right)\bigg)_{n \in \nn} =
\bigg(\bigcap_{\lambda \in \Lambda}\a_{\lambda,n}\bigg)_{n \in
\nn}. $$ The final claim now follows from \ref{pd.1j}.
\end{proof}

We can deduce quickly from Theorem \ref{pd.2} that $R^{\infty}$ is
only Noetherian in rather uninteresting cases.

\begin{lem}
\label{pd.3}
Let $\fP \in \Spec (R^{\infty})$. Then $\fP^{[p^n]} = \fP^n = \fP$ for
all $n \in \N$.
\end{lem}

\begin{proof} Let $\fp \in \Spec (R)$ correspond to $\fP$ in the
correspondence of \ref{pd.2}(i). Let $\alpha \in \fP$. With the
notation of \ref{pd.1j}, we have $\alpha = x^{-i}\overline{r}x^i$
for some $r \in \fp$ and some $i \in \nn$. Now
$x^{-i-1}\overline{r}x^{i+1} \in x^{-i-1}\pi(\fp) x^{i+1}
\subseteq \fP$. Hence
$
\alpha = x^{-i}\overline{r}x^i = x^{-i-1}\overline{r}^px^{i+1} =
(x^{-i-1}\overline{r}x^{i+1})^p \in \fP^{[p]}.
$
Hence $\fP \subseteq \fP^{[p]}$. It follows that $\fP =
\fP^{[p]}$, and then, by induction, that $\fP = \fP^{[p^n]}$ for
all $n \in \N$.

Finally, $\fP \supseteq \fP^n \supseteq \fP^{[p^n]} = \fP$ for all
$n \in \N$.
\end{proof}

We next justify our claim that $R^{\infty}$ is only Noetherian in
uninteresting cases.

\begin{thm}
\label{pd.4}
The following conditions are equivalent:

\begin{enumerate}
\item $R^\infty$ is Noetherian;
\item $\R$ is a direct product of a finite number of fields;
\item $R^\infty$ is a direct product of a finite number of fields.
\end{enumerate}
\end{thm}

\begin{proof} We can, and do, assume that $R$ is reduced.

(i) $\Rightarrow$ (ii) As $R^\infty$ is Noetherian,
  every prime ideal of $R^\infty$ is finitely generated.
  As for every prime ideal $\fP$ of $R^\infty\/,$ we have $\fP^2 =
  \fP$ (by \ref{pd.3}),
  Nakayama's Lemma implies that each prime ideal in $R^\infty$ is
  a minimal prime. Hence $\dim {R^\infty} = 0\/.$
  Therefore $R^\infty$ is Artinian with finitely many maximal ideals,
  $\fM_1, \dots, \fM_n$ say. By Theorem \ref{pd.2}(i),
  $\Spec R = \{ \fm_1, \dots, \fm_n \}$
  where $\fm_i$ corresponds to $\fM_i$ for $i = 1, \ldots, n$ and
  each $\fm_i$ ($1 \leq i \leq n$) is a maximal ideal of $R\/.$
  As $R$ is reduced, the natural ring homomorphism
  $R \lra R/\fm_1 \times
  \dots \times R/\fm_n$ is an isomorphism.

  (ii) $\Rightarrow$ (iii) This follows from \ref{pc.3} and
  \cite[Proposition 3.1]{Jorda82}.

  (iii) $\Rightarrow$ (i) This is clear.
  \end{proof}

Thus it is not clear whether a given proper ideal $\fA$ of
$R^{\infty}$ has a primary decomposition (although we shall
present in \S\ref{ex} some examples of proper ideals in
non-Noetherian perfect closures that \emph{do} have primary
decompositions). But if it does, we can draw some interesting
conclusions, as we now show.

\begin{thm}
\label{pd.5} Let $\fA$ be an ideal of $R^\infty$ and let
$(\fa_n)_{n \in \nn}$
 be the corresponding $f$-sequence of ideals of $R$. For
each $k = 1, \ldots, t$, let $\fP_k$ be a prime ideal of
$R^\infty$, let $\fp_k$ be the corresponding prime ideal of $R$
(see\/ {\rm \ref{pd.2}(i)}), let $\fQ_k$ be an ideal of
$R^{\infty}$, and let $(\fq_{k,n})_{n \in \nn}$ be the
$f$-sequence of ideals of $R$ corresponding to $\fQ_k$.

 Then the following are equivalent:
 \begin{enumerate}
  \item $\fA$ has a primary decomposition
    $\fA = \fQ_1 \cap \dots \cap \fQ_t$ where $\fQ_k$ is $\fP_k$-primary
    for $k = 1, \ldots, t$;
  \item each $\fa_n$ has a primary decomposition
    $\fa_n = \fq_{1,n} \cap \dots \cap \fq_{t,n}$ such that
    $\fq_{k,n}$ is $\fp_k$-primary
       for all $k = 1, \ldots, t$ and all $n \in \nn$.
 \end{enumerate}
Furthermore, when these equivalent conditions are satisfied,
$$\bigcup_{n\in \nn} \ass \fa_n = \left\{ \pi^{-1}(\fP \cap \R) : \fP \in
\ass_{R^{\infty}}(\fA) \right\} = \ass \fa_j \quad \mbox{~for
all~} j \gg 0,$$ and so is finite; furthermore, the primary
decomposition of $\fA$ given in\/ {\rm (i)} is minimal if and only
if the primary decomposition of $\fa_n$ given in\/ {\rm (ii)} is
minimal for all $n \gg 0$.
\end{thm}

\begin{proof} Assume that statement (i) is true.
By \ref{pd.2}(iv), $\fq_{k,n}$ is $\fp_k$-primary
       for all $k = 1, \ldots, t$ and all $n \in \nn$; furthermore,
       with the notation of \ref{pd.1j}, we have, for each $n \in \nn$,
\[
\fa_n = \psi_n^{-1}(\fA) = \psi_n^{-1}\bigg( \bigcap_{k=1}^t \fQ_k
\bigg) = \bigcap_{k=1}^t\psi_n^{-1}\left( \fQ_k \right) =
\bigcap_{k=1}^t\fq_{k,n}.
\]
Thus statement (ii) is true.

Now suppose, in addition, that $\fA = \fQ_1 \cap \dots \cap \fQ_t$ is a
minimal primary decomposition of $\fA$. The bijective correspondence of
\ref{pd.2}(i) shows that $\fp_1, \ldots, \fp_t$ are all distinct. Choose
$j \in \{1, \ldots,t \}$. Since (in view of \ref{pd.2}(v))
\[
\bigcup_{n \in \nn} x^{-n} \pi(\fq_{j,n})x^n = \fQ_{j}
\not\supseteq \bigcap_{\stackrel{\scriptstyle k=1}{k \neq j}}^{t}
\fQ_{k} = \bigcup_{n \in \nn} \left(x^{-n}\pi\left(
\bigcap_{\stackrel{\scriptstyle k=1}{k \neq j}}^{t} \fq_{k,n}
\right) x^n\right),
\]
there exists $n \in \nn$ such that $$\fq_{j,n} \not\supseteq
\bigcap_{\stackrel{\scriptstyle k=1}{k \neq j}}^{t} \fq_{k,n}.$$
Hence refinement of the primary decomposition $\fa_n = \fq_{1,n}
\cap \dots \cap \fq_{t,n}$ to a minimal primary decomposition
cannot result in the removal of $\fq_{j,n}$. Therefore $\fp_j \in
\ass \a_n$. Since the sequence $\left(\ass \a_n\right)_{n \in
\nn}$ is increasing (by \ref{pd.1a}(iv)), it follows that $\fp_j \in
\ass \a_{n+h}$ for all $h \in \nn$. It follows that the primary
decomposition of $\fa_m$ given in (ii) is minimal for all $m \gg
0$.

Now assume that statement (ii) is true. It follows from
\ref{pd.2}(iv) that $\fQ_k$ is $\fP_k$-primary, and from
\ref{pd.2}(v) that
$$
\fA = \bigcup_{n \in \nn} x^{-n}\pi(\fa_n)x^n = \bigcup_{n \in \nn} \bigg(x^{-n}
\pi\bigg(\bigcap_{k=1}^t\fq_{k,n}\bigg)x^n\bigg) = \bigcap_{k=1}^t \bigg(
\bigcup_{n \in \nn} x^{-n}\pi(\fq_{k,n})x^n\bigg) = \bigcap_{k=1}^t \fQ_k.
$$
Thus statement (i) is true.

Now suppose that the primary decomposition of $\fa_j$ given
in (ii) is minimal for all $j \gg 0$. This means that
$\fp_1, \ldots, \fp_t$ are all different (so that
$\fP_1, \ldots, \fP_t$ are all different) and $\ass \fa_j =
\{\fp_1, \ldots, \fp_t\}$ for all $j \gg 0$. Now if the primary decomposition
$\fA = \bigcap_{k=1}^t \fQ_k$ were not minimal, then
it would be possible to refine it to a minimal primary decomposition which
would have fewer than $t$ terms; it would then follow from what we have proved
above (when dealing with the implication (i) $ \Rightarrow $ (ii)) that
all the sets $\ass \a_n~(n \in \nn)$ would have fewer than $t$ terms, and this
would be a contradiction. Hence $\fA = \bigcap_{k=1}^t \fQ_k$ must be a
minimal primary decomposition.
\end{proof}

It follows from Theorem \ref{pd.4} that, provided $\R$ is not a
direct product of a finite number of fields, there will be ideals
in $R^{\infty}$ that are not finitely generated. It will be
helpful to have information about exactly when an ideal of
$R^{\infty}$ is finitely generated, and our next result provides
this. Its proof uses some ideas from the proof of Jordan
\cite[Theorem 5.6]{Jorda82}. Observe that if $\fA = (\alpha_1,
\dots, \alpha_t)R^\infty$ is a finitely generated ideal of
$R^\infty\/,$ then there exist $k \in \nn$ and $a_1, \dots, a_t
\in R$ such that $\alpha_j = x^{-k}\overline{a_j}x^k$ for all $j = 1, \ldots,
t\/.$

\begin{thm}
\label{ex.1}
 Let $\fA \subseteq R^\infty$ be an ideal and $(\fa_n)_{n \in\nn}$
 be the $f$-sequence corresponding to $\fA\/.$ Let $k \in \nn$ and
$a_1, \ldots, a_t \in R$.
 Then, with the notation of \/ {\rm \ref{pd.1j}},
 the following statements are equivalent:
 \begin{enumerate}
  \item $\fA = (x^{-k}\overline{a_1}x^k, \dots, x^{-k}\overline{a_t}x^k)R^\infty$;
  \item $\fa_{k+n} = ((a_1, \dots, a_t)^{[p^n]})^F$
  for all $n \in \nn$.
 \end{enumerate}
\end{thm}

\begin{note} When statements (i) and (ii) in the theorem are satisfied,
we can describe in precise terms not only the $\fa_{k+n}$ for $n \in \nn$,
but also each $\fa_j$ for $0 \leq j < k$: since $(\fa_n)_{n \in\nn}$ is
an $f$-sequence, we must have $\fa_j = f^{-(k-j)}(\fa_k) = f^{-(k-j)}
\big((a_1, \dots, a_t)^F\big)$.

Note also that a special case of this theorem yields that, if $\fa$ is
an ideal of $R$, then the $f$-sequence that corresponds to the ideal
$\fa R^{\infty}$ of $R^{\infty}$ is $\big((\fa^{[p^n]})^F\big)_{n \in \nn}$.
\end{note}

\begin{proof} (i) $\Rightarrow$ (ii) Let $n \in \nn$.
 As $x^{-k}\overline{a_j}x^k \in \fA$, we have
 $a_j \in \fa_k$ for all $ j = 1, \ldots,
t\/,$ and so it follows from
 \ref{pd.1a}(i) that
 $
  ((a_1, \dots, a_t)^{[p^n]})^F
  \subseteq \fa_{k+n}\/.
 $

 Next consider $b \in \fa_{k+n}\/.$
 Then, by definition, $x^{-(k+n)}\overline{b}x^{k+n} \in \fA\/.$
 As $\fA = (x^{-k}\overline{a_1}x^k, \dots, x^{-k}\overline{a_t}x^k)R^{\infty}$,
 we can write $x^{-(k+n)}\overline{b}x^{k+n} = \sum_{j = 1}^t
 \rho_jx^{-k}\overline{a_j}x^k$
 for some $\rho_1, \ldots,\rho_t \in R^\infty\/.$
 Hence
 $$
  \overline{b} = \sum_{j = 1}^t \rho_j^{p^{k+n}}x^{-k}\overline{a_j}^{p^{k+n}}x^k
    = \sum_{j = 1}^t \rho_j^{p^{k+n}}\overline{a_j}^{p^n}\/.
 $$
 As $\rho_1, \dots, \rho_t \in R^\infty\/,$ there exists $m \in
 \nn$ such that, for all $j = 1, \ldots, t$, we have
 $\rho_j^{p^{k+n+m}} \in \R\/.$
 Therefore $\overline{b}^{p^m} = \sum_{j = 1}^t
 \rho_j^{p^{k+n+m}}\overline{a_j}^{p^{n+m}} \in
 (\overline{a_1}, \dots,  \overline{a_t})^{[p^{n+m}]}\/.$
 Hence there exists $b' \in (a_1, \ldots, a_t)^{[p^{n+m}]}$ such that $b^{p^m}
 - b' \in \sqrt{(0)}$, and it follows easily that $b
 \in ((a_1, \dots, a_t)^{[p^n]})^F\/.$
 Therefore
 $$
   \fa_{k+n} \subseteq ((a_1, \dots,
    a_t)^{[p^n]})^F\/.
 $$

 (ii) $\Rightarrow$ (i)
 As $\fa_k = (a_1, \dots, a_t)^F$ and $\fA = \bigcup_{n \in \nn}
 x^{-n}\pi(\fa_n)x^n\/,$ it follows that
 $x^{-k} \overline{a_1}x^k, \dots, x^{-k} \overline{a_t}x^k
 \in \fA\/.$
 Hence $(x^{-k} \overline{a_1}x^k, \dots, x^{-k} \overline{a_t}x^k)R^\infty
 \subseteq \fA\/.$

 Next consider $\gamma \in \fA$.
 We can write $\gamma = x^{-i} \overline{a}x^i$ for some $i \in \nn$ and some
 $a \in \fa_i\/.$ If $i < k\/,$ then $\gamma =
 x^{-k} \overline{a^{p^{k-i}}}x^k$
 with $a^{p^{k-i}} \in \fa_k\/.$ Hence we can assume that $i \geq
 k\/$. Write $i = k+n$, where $n \in \nn$.

 Then, as $a \in \fa_{k+n} = ((a_1, \dots, a_t)^{[p^n]})^F\/,$
 there exists $m \in \nn$ such that $$a^{p^m} \in ((a_1, \dots,
 a_t)^{[p^n]})^{[p^m]} = (a_1, \dots, a_t)^{[p^{n+m}]}\/.$$
 Therefore, we can write
 $
  a^{p^m} = \sum_{j = 1}^t r_j a_j^{p^{n+m}}
 $
 for some $r_1, \dots, r_t \in R\/.$
 Hence
 \begin{align*}
   \gamma^{p^m} &= (x^{-(k+n)} \overline{a}x^{k+n})^{p^m}
   = x^{-(k+n)} \overline{a}^{p^m}x^{k+n}
   = x^{-(k+n)}\bigg(\sum_{j = 1}^t  \overline{r_j
    a_j^{p^{n+m}}}\bigg)x^{k+n}\\
   &= \sum_{j = 1}^t x^{-(k+n)} \overline{r_j}x^{k+n}
   x^{-(k+n)} \overline{a_j}^{p^{n+m}}x^{k+n}
   = \sum_{j = 1}^t x^{-(k+n+m)} \overline{r_j}^{p^m}x^{k+n+m}
       x^{-k} \overline{a_j}^{p^{m}}x^k\\
   &= \bigg(\sum_{j = 1}^t x^{-(k+n+m)} \overline{r_j}x^{k+n+m}
       x^{-k} \overline{a_j}x^k\bigg)^{p^m}.
 \end{align*}
 Therefore $
  \gamma = \sum_{j = 1}^t x^{-(k+n+m)} \overline{r_j}x^{k+n+m}
       x^{-k} \overline{a_j}x^k
 $
 because the Frobenius homomorphism acts bijectively on $R^\infty\/$;
hence $\gamma \in (x^{-k} \overline{a_1}x^k, \dots,
 x^{-k} \overline{a_t}x^k)R^\infty$.
\end{proof}

\section{\sc Linear growth of primary decompositions}
\label{lg}

Throughout this section, $R$ will denote a commutative
Noetherian ring of prime characteristic $p$. In the following
definition, we extend the concept of `linear growth of primary
decompositions', defined in the Introduction for the sequence of
Frobenius powers of a given ideal of $R$, to $f$-sequences.

\begin{defi}
\label{lg.1} Let $(\fa_n)_{n \in\nn}$ be an $f$-sequence of proper
ideals of $R$ and let $h$ be a positive integer.

We say that $(\fa_n)_{n \in\nn}$ has {\em $h$-linear growth of
primary decompositions\/} if, for every non-negative integer $n$,
there exists a minimal primary decomposition $\fa_n = \fq_{1,n}
\cap \ldots \cap \fq_{k_n,n}$ with $\sqrt{\fq_{i,n}}^{[p^n]h}
\subseteq \fq_{i,n}$ for all $i =1, \ldots, k_n$. We say that
$(\fa_n)_{n \in\nn}$ has {\em linear growth of primary
decompositions\/} if it has $k$-linear growth of primary
decompositions for some positive integer $k$.

Let $\fa$ be a proper ideal of $R$. Although the sequence of
Frobenius powers of $\fa$ need not be an $f$-sequence, we say that
the Frobenius powers of $\fa$ have {\em $h$-linear growth of
primary decompositions\/} if, for every non-negative integer $n$,
there exists a minimal primary decomposition $$\fa^{[p^n]} = \fq_{1,n}
\cap \ldots \cap \fq_{k_n,n} \quad \mbox{~with~} \sqrt{\fq_{i,n}}^{[p^n]h}
\subseteq \fq_{i,n} \mbox{~for all~} i =1, \ldots, k_n.$$
\end{defi}

An elementary argument will prove the following lemma.

\begin{lem}
\label{bk.44} Let $\eta : R \lra R'$ be a homomorphism of
commutative Noetherian rings of characteristic $p$, and let
$(\fa'_n)_{n \in\nn}$ be an $f$-sequence of ideals of $R'$. Let
$\eta^{\infty} : R^{\infty} \lra R'^{\infty}$ be the ring
homomorphism induced by $\eta$ (see\/ {\rm \S1}). By\/ {\rm
\ref{bk.4}}, the sequence $\left(\eta^{-1}(\fa'_n)\right)_{n
\in\nn}$ is an $f$-sequence of ideals of $R$.

If $(\fa'_n)_{n \in\nn}$ has $h$-linear growth of
primary decompositions for some $h \in \N$, then
$\left(\eta^{-1}(\fa'_n)\right)_{n \in\nn}$ also
has $h$-linear growth of
primary decompositions.
\end{lem}

Several of the results in this section will involve the hypothesis
that $R$ has a $p^{m_0}$-weak test element for some $m_0 \in \nn$.
We remind the reader of the definition of such elements.

\begin{defi}
\label{lg.32} (See \cite[\S6]{HocHun90}.) Let $m_0 \in \nn$. A
{\em $p^{m_0}$-weak test element for $R$\/} is an element $c \in
R^{\circ}$ such that, for every ideal $\fb$ of $R$ and for $r \in
R$, it is the case that $r \in \fb^*$ if and only if $cr^{p^n} \in
\fb^{[p^n]}$ for all $n \geq m_0$. A $p^0$-weak test element for
$R$ is called a {\em test element for $R$\/}.

Such a $p^{m_0}$-weak test element $c$ for $R$ is said to be a
{\em locally stable $p^{m_0}$-weak test element\/} for $R$ if, for
each $\fp \in \Spec (R)$, the image of $c$ in the localization
$R_{\fp}$ is a  $p^{m_0}$-weak test element for $R_{\fp}$.

Also, a $p^{m_0}$-weak test element $c$ for $R$ is said to be a
{\em completely stable $p^{m_0}$-weak test element\/} for $R$ if
it is locally stable and, for each $\fp \in \Spec (R)$, the image
of $c$ in the completion $\widehat{R_{\fp}}$ of the localization
of $R$ at $\fp$ is a $p^{m_0}$-weak test element for
$\widehat{R_{\fp}}$.

It is a result of Hochster and Huneke \cite[Theorem
(6.1)(b)]{HocHun94} that an algebra of finite type over an
excellent local ring of characteristic $p$ has a completely stable
$p^{m_0}$-weak test element, for some $m_0 \in \nn$.
\end{defi}

\begin{lem}
\label{bk.3} Let $\eta : R \lra R'$ be a homomorphism of
commutative Noetherian rings of characteristic $p$, let $S$ be a
multiplicatively closed subset of $R$, and let $\phantom{\fa}^{e}$
stand for extension with respect to the natural ring homomorphism
$R \lra S^{-1}R$.

\begin{enumerate}
\item If $\fa'$ is an $F$-closed ideal of $R'$, then $\eta^{-1}(\fa')$
is an $F$-closed ideal of $R$.
\item If $\fa'$ is a tightly-closed ideal of $R'$ and there exists a
$p^{m_0}$-weak test element (where $m_0 \in \nn$) $c$ for $R$ such
that $\eta(c)$ is a $p^{m_0}$-weak test element for $R'$, then
$\eta^{-1}(\fa')$ is a tightly-closed ideal of $R$.
\item If $\fa$ is an $F$-closed ideal of $R$, then $\fa^e$
is an $F$-closed ideal of $S^{-1}R$.
\end{enumerate}
\end{lem}

\begin{proof} (i) Let $r \in (\eta^{-1}(\fa'))^F$, so that there exists $n \in \nn$
such that $r^{p^n} \in (\eta^{-1}(\fa'))^{[p^n]}$. Thus there
exist $a_1, \ldots, a_t \in \eta^{-1}(\fa')$ and $r_1, \ldots, r_t
\in R$ such that
$
r^{p^n} = r_1a_1^{p^n} + \cdots +  r_ta_t^{p^n}.
$
Apply $\eta$ to see that $(\eta (r))^{p^n} \in \fa'^{[p^n]}$, so
that $\eta (r) \in \fa'$ because $\fa'$ is $F$-closed.

(ii) Let $r \in (\eta^{-1}(\fa'))^*$, so that $cr^{p^n} \in
(\eta^{-1}(\fa'))^{[p^n]}$ for all $n \geq m_0$. Application of
$\eta$ then yields that $\eta(c)\eta(r)^{p^n} \in \fa'^{[p^n]}$
for all $n \geq m_0$. Since $\eta(c)$ is a $p^{m_0}$-weak test
element for $R'$ and $\fa'$ is tightly closed, $\eta(r) \in \fa'$.

(iii) Let $r \in R$ be such that $r/1 \in (\fa^e)^F$ in $S^{-1}R$,
so that there exists $n \in \nn$ with $r^{p^n}/1 \in
(\fa^e)^{[p^n]} = (\fa^{[p^n]})^e$. Thus there exists $s \in S$
such that $sr^{p^n} \in \fa^{[p^n]}$; since $sr \in \fa^F = \fa$,
it follows that $r/1 = sr/s \in \fa^e$.
\end{proof}

The motivation for our study of linear growth of primary
decompositions came from work of K. E. Smith and I. Swanson: in
\cite[Theorem (1.2)]{KESISLinGrowth}, they showed that, in $K[X_1,
\ldots,X_d]/\fb$, where $K$ is a field of characteristic $p$,
$X_1, \ldots,X_d$ are indeterminates, and $\fb$ is an ideal of the
polynomial ring $K[X_1, \ldots,X_d]$ generated by monomials, the
Frobenius powers of a proper ideal $\fa$ generated by monomials in
the images of $X_1, \ldots,X_d$ have linear growth of primary
decompositions (in the sense of Definition \ref{lg.1}), and\/
$\bigcup_{n \in \nn} \ass \fa^{[p^n]}$ is a finite set.

Our first major result in this section
demonstrates a connection between the two definitions in
\ref{lg.1}: it shows that, if the Frobenius powers of the
proper ideal $\fa$ of $R$ have $h$-linear growth of primary
decompositions and $\bigcup_{n \in \nn} \ass \fa^{[p^n]}$ is
finite, then the canonical $f$-sequence
$\big((\fa^{[p^n]})^F\big)_{n \in\nn}$ associated to $\fa$ (see
\ref{pd.0}(iii)) has $h$-linear growth of primary
decompositions. However, this result should come with the warning
that M. Katzman \cite{Katzm96} has provided an example of a proper
ideal $\fd$ in a ring of the type under consideration here for
which $\bigcup_{n \in \nn} \ass \fd^{[p^n]}$ is infinite.

\begin{thm}
\label{lg.2} Suppose that the Frobenius powers of the proper ideal
$\fa$ of $R$ have $h$-linear growth of primary decompositions (for
some $h \in \N$), and that\/ $\bigcup_{n \in \nn} \ass
\fa^{[p^n]}$ is finite; let the members of the latter set be
$\fp_1, \ldots, \fp_t$, and let the corresponding prime ideals of
$R^{\infty}$ (see\/ {\rm \ref{pd.2}(i)}) be $\fP_1, \ldots, \fP_t$
respectively.

For each $i = 1, \ldots, t$, let $\phantom{\fa}^{e_i}$ and
$\phantom{\fa}^{c_i}$ stand for extension and
 contraction with respect to the natural ring homomorphism
$R \lra R_{\fp_i}$, and let $\phantom{\fa}^{E_i}$ and
$\phantom{\fa}^{C_i}$ stand for extension and
 contraction with respect to the natural ring homomorphism
$R^{\infty} \lra \left(R^{\infty}\right)_{\fP_i}$.

\begin{enumerate}
\item For each $n \in \nn$,
$$
\fa^{[p^n]} = \bigcap_{i=1}^t \big((\fa +
\fp_i^h)^{[p^n]}\big)^{e_ic_i}, \quad \mbox{~where~} \big((\fa +
\fp_i^h)^{[p^n]}\big)^{e_ic_i} \mbox{~is~} \fp_i\mbox{-primary
for~} i = 1, \ldots,t,
$$
is a primary decomposition.
\item For each $n \in \nn$,
$$
(\fa^{[p^n]})^F = \bigcap_{i=1}^t \big(((\fa +
\fp_i^h)^{[p^n]})^F\big)^{e_ic_i}, \quad \mbox{~where~} \big(((\fa
+ \fp_i^h)^{[p^n]})^F\big)^{e_ic_i} \mbox{~is~}
\fp_i\mbox{-primary for~} i = 1, \ldots,t,
$$
is a primary decomposition in which all the primary components are
$F$-closed.
\item The $f$-sequence $\big((\fa^{[p^n]})^F\big)_{n \in \nn}$
has $h$-linear growth of primary decompositions.
\item We have $\bigcup_{n \in \nn} \ass (\fa^{[p^n]})^F \subseteq
\bigcup_{n \in \nn} \ass \fa^{[p^n]}$, and so $\bigcup_{n \in \nn}
\ass (\fa^{[p^n]})^F$ is finite.
\item Furthermore,
$$
\fa R^{\infty} = \bigcap_{i=1}^t \big((\fa +
\fp_i^h)R^{\infty}\big)^{E_iC_i}, \quad \mbox{~where~} \big((\fa +
\fp_i^h)R^{\infty}\big)^{E_iC_i} \mbox{~is~} \fP_i\mbox{-primary
for~} i = 1, \ldots,t,
$$
is a primary decomposition of the ideal $\fa R^{\infty}$ of $R^{\infty}$ to which
the $f$-sequence $\big((\fa^{[p^n]})^F\big)_{n \in \nn}$ corresponds.
\end{enumerate}
\end{thm}

\begin{proof} For each $n \in \nn$, there exists a subset
$\Lambda_n$ of $\{1, \ldots, t\}$ and a minimal primary
decomposition $\fa^{[p^n]} = \bigcap_{i \in \Lambda_n} \fq_{i,n}$
where $\fq_{i,n}$ is $\fp_i$-primary and $\fp_i^{h[p^n]} \subseteq
\fq_{i,n}$ for all $i \in \Lambda_n$.

(i) Now $\fa^{[p^n]} \subseteq \fa^{[p^n]} + \fp_i^{h[p^n]} = (\fa
+ \fp_i^h)^{[p^n]} \subseteq \fq_{i,n},$ so that $$\fa^{[p^n]}
\subseteq \big((\fa + \fp_i^h)^{[p^n]}\big)^{e_ic_i} \subseteq
(\fq_{i,n})^{e_ic_i} = \fq_{i,n} \quad \mbox{~for all~} i \in
\Lambda_n.$$ Hence $\fa^{[p^n]} = \bigcap_{i \in
\Lambda_n}\big((\fa + \fp_i^h)^{[p^n]}\big)^{e_ic_i}$; moreover,
for all $i \in \{1, \ldots,t\}$, we have $\fp_i^{h[p^n]} \subseteq
(\fa + \fp_i^h)^{[p^n]} \subseteq \fp_i$, so that $\big((\fa +
\fp_i^h)^{[p^n]}\big)^{e_ic_i}$ is $\fp_i$-primary. It follows
that
$$
\fa^{[p^n]} = \bigcap_{i=1}^t \big((\fa +
\fp_i^h)^{[p^n]}\big)^{e_ic_i}, \quad \mbox{~where~} \big((\fa +
\fp_i^h)^{[p^n]}\big)^{e_ic_i} \mbox{~is~} \fp_i\mbox{-primary
for~} i = 1, \ldots,t,
$$
is a primary decomposition.

(ii),(iii),(iv) It is clear that $(\fa^{[p^n]})^F \subseteq
\bigcap_{i=1}^t ((\fa + \fp_i^h)^{[p^n]})^F \subseteq
\bigcap_{i=1}^t \big(((\fa + \fp_i^h)^{[p^n]})^F \big)^{e_ic_i}$.
Note that, for each $i \in \{1, \ldots,t\}$, we have $
\fp_i^{h[p^n]} \subseteq (\fa + \fp_i^h)^{[p^n]}  \subseteq
\fp_i$, so that
$$
\fp_i^{h[p^n]} \subseteq ((\fa + \fp_i^h)^{[p^n]})^F \subseteq
\big(((\fa + \fp_i^h)^{[p^n]})^F \big)^{e_ic_i} \subseteq
\big(\fp_i^F\big)^{e_ic_i} = \fp_i^{e_ic_i} = \fp_i.
$$
It therefore follows that $\big(((\fa + \fp_i^h)^{[p^n]})^F
\big)^{e_ic_i}$ is $\fp_i$-primary; it is $F$-closed by
\ref{bk.3}(i),(iii). Therefore, the proof of parts
(ii), (iii) and (iv) will have been
completed as soon as it has been
shown that $ (\fa^{[p^n]})^F \supseteq \bigcap_{i=1}^t
\big(((\fa + \fp_i^h)^{[p^n]})^F \big)^{e_ic_i} $ for each $n \in
\nn$.

So let $r \in \bigcap_{i=1}^t \big(((\fa + \fp_i^h)^{[p^n]})^F
\big)^{e_ic_i}.$ Then, for each $i \in \{1, \ldots, t\}$, there
exists $s_i \in R \setminus \fp_i$ such that $s_ir \in ((\fa +
\fp_i^h)^{[p^n]})^F$; this means that there exists $k_i \in \nn$
with $(s_ir)^{p^{k_i}} \in ((\fa + \fp_i^h)^{[p^n]})^{[p^{k_i}]}$.
Set $k := \max \{k_i:i = 1, \ldots, t\}$. Then $(s_ir)^{p^{k}} \in
((\fa + \fp_i^h)^{[p^n]})^{[p^{k}]} = (\fa +
\fp_i^h)^{[p^{n+k}]}$, so that
$$
r^{p^{k}} \in \big((\fa + \fp_i^h)^{[p^{n+k}]}\big)^{e_ic_i} \quad
\mbox{~for all~} i = 1, \ldots, t.
$$
Hence, on use of part (i), we see that
$
r^{p^{k}} \in \bigcap_{i=1}^t\big((\fa +
\fp_i^h)^{[p^{n+k}]}\big)^{e_ic_i} = \fa^{[p^{n+k}]} =
(\fa^{[p^{n}]})^{[p^{k}]}.
$
Therefore $r \in (\fa^{[p^n]})^F$, as required.

(v) For each $i = 1, \ldots, t$, the sequence $\big(((\fa +
\fp_i^h)^{[p^n]})^F\big)_{n \in \nn}$ is the canonical
$f$-sequence associated to $\fa + \fp_i^h$: see \ref{pd.0}(iii).
This $f$-sequence corresponds to $(\fa + \fp_i^h)R^{\infty}$, by
Theorem \ref{ex.1}; similarly, $\fa R^{\infty}$ is the ideal of
$R^{\infty}$ corresponding to the $f$-sequence
$\big((\fa^{[p^n]})^F\big)_{n \in \nn}$. It is easy to check that $\bigcup_{n \in
\nn} \psi_n(R \setminus \fp_i) = R^{\infty} \setminus \fP_i$,
and so it follows from Lemma \ref{pd.2a} that $((\fa +
\fp_i^h)R^{\infty})^{E_iC_i}$ corresponds to the $f$-sequence
$\big(\big(((\fa + \fp_i^h)^{[p^n]})^F\big)^{e_ic_i}\big)_{n \in
\nn}$. By part (ii), each term in the latter $f$-sequence is
$\fp_i$-primary. The result therefore follows from
Theorem \ref{pd.5}.
\end{proof}

It was explained in the Introduction that our interest in linear
growth of primary decompositions of Frobenius powers of a proper
ideal $\fa$ of $R$ arose from the argument of Smith and Swanson in
\cite{KESISLinGrowth} that shows that if the Frobenius powers of
$\fa$ have linear growth of primary decompositions, then $\fa^*R_u
= (\fa R_u)^*$ for all $u \in R$. Our next result shows, among
other things, that linear growth of primary decompositions of an
$f$-sequence $(\fa_n)_{n \in\nn}$ that approximates to the
Frobenius powers of $\fa$ in the sense that $\fa^{[p^n]} \subseteq
\fa_n \subseteq(\fa^{[p^n]})^*$ for all $n \in \nn$ would do just
as well in this context, provided that $R$ has a $p^{m_0}$-weak
test element, for some $m_0 \in \nn$. Note that Theorem \ref{lg.2}
shows that if the Frobenius powers of $\fa$ have linear growth of
primary decompositions, then, provided $\bigcup_{n \in \nn}\ass
\fa^{[p^n]}$ is a finite set, there {\em is\/} such an
$f$-sequence with linear growth of primary decompositions, namely
$((\fa^{[p^n]})^F)_{n \in\nn}$.

\begin{thm}
\label{lg.3} Assume that $R$ has a  $p^{m_0}$-weak test element for
some $m_0 \in \nn$; let $\fa$ be a proper
ideal of $R$. Suppose that $(\fa_n)_{n \in\nn}$ is an $f$-sequence
of ideals of $R$ which has linear growth of primary decompositions
and is such that $\fa^{[p^n]} \subseteq \fa_n
\subseteq(\fa^{[p^n]})^*$ for all $n \in \nn$.
\begin{enumerate}
\item For each $u \in R$, we have $\fa^*R_u = (\fa R_u)^*$.
\item If\/ $\bigcup_{n \in \nn}\ass \fa_n$ is a finite set, then $\fa^*S^{-1}R =
(\fa S^{-1}R)^*$ for every multiplicatively closed subset $S$ of
$R$.
\end{enumerate}
\end{thm}

\begin{note} Note that, if $(\fa'_n)_{n \in\nn}$ is
an $f$-sequence of ideals of $R$ such that $\fa \subseteq \fa'_0$,
then, for each $n \in \nn$, we have $\fa^{[p^n]} \subseteq
(\fa'_0)^{[p^n]} \subseteq \fa'_n$, by \ref{pd.1a}(i).

Also, the condition that $\fa^{[p^n]} \subseteq \fa_n \subseteq(\fa^{[p^n]})^*$
(for some $n \in \nn$) can be described, in
the terminology of \cite[(7.11)]{HocHun90}, by saying that
$\fa_n$ is {\it trapped\/} over $\fa^{[p^n]}$.
\end{note}

\begin{proof} Let $S$ denote an
arbitrary multiplicatively closed subset of $R$.
It is clear that $\fa^*S^{-1}R \subseteq (\fa S^{-1}R)^*$.

Let $c' \in R^{\circ}$ be a
$p^{m_0}$-weak test element for $R$.  Let
$r \in R$ be such that $r/1 \in (\fa
S^{-1}R)^*$. Thus there exists $c \in R$ such that $c/1 \in
\left(S^{-1}R\right)^\circ$ and $cr^{p^n}/1 \in (\fa
S^{-1}R)^{[p^n]}$ (in the ring of fractions $S^{-1}R$) for all $n
\gg 0$. If $c \not\in R^\circ$, let $d$ be an element of $R$
that belongs to those minimal primes of $R$ to which $c$ does not
belong, and to no others: an argument in the proof of
\cite[Proposition 4.14]{HocHun90} shows that we may add a suitable
power of $d$ to $c$ to see that we can assume that $c \in
R^\circ$.

By hypothesis, there exists $h \in\N$ such that,
for each $n \in \nn$, there is a minimal primary decomposition
$\fa_n = \fq_{1,n} \cap \ldots \cap
\fq_{k_n,n}$ with $\sqrt{\fq_{i,n}}^{[p^n]h} \subseteq \fq_{i,n}$
for all $i =1, \ldots, k_n$.

(i) Here we consider the special case in which $S = \{u^k:k \in \nn\}$, and
our conclusion above specialises to the statement that $cr^{p^n}/1 \in (\fa
R_u)^{[p^n]}$ (in the ring of fractions $R_u$) for all $n
\gg 0$, say for all $n \geq n_0$.

Therefore, for each $n \geq n_0$, there exists $t(n) \in \nn$ such
that $cu^{t(n)}r^{p^n} \in \fa^{[p^n]}$. Choose $n \geq n_0$ and
$i \in \{1, \ldots, k_n\}$. Now $cu^{t(n)}r^{p^n} \in \fa^{[p^n]}
\subseteq \fa_n \subseteq \fq_{i,n}$. If $cr^{p^n} \not\in
\fq_{i,n}$, then $u^{t(n)} \in \sqrt{\fq_{i,n}}$, so that $u \in
\sqrt{\fq_{i,n}}$ and $u^{p^nh} \in \fq_{i,n}$. It follows that $
c(ru^h)^{p^n} \in  \bigcap_{i=1}^{k_n} \fq_{i,n} = \fa_n
\subseteq(\fa^{[p^n]})^*. $ It should be noted that this is true
for each $n \geq n_0$. Since $c'$ is a $p^{m_0}$-weak test element
for $R$, we have $c'(c(ru^h)^{p^n})^{p^{m_0}} \in
(\fa^{[p^n]})^{[p^{m_0}]}$ for all $n  \geq n_0$. Therefore $
c'c^{p^{m_0}}(ru^h)^{p^{n+m_0}} \in \fa^{[p^{n+m_0}]}$ for all $n
\geq n_0.$ Since $c'c^{p^{m_0}} \in R^{\circ}$, it follows that
$ru^h \in \fa^*$, so that $r/1 = ru^h/u^{h} \in \fa^*R_u$.

(ii) Here, we revert to the situation where $S$ is an arbitrary
multiplicatively closed
subset of $R$.

Note that, by \ref{pd.1a}(iv), $\ass \fa_n \subseteq \ass
\fa_{n+1}$ for all $n \in \nn$. Suppose that the finite set
$\bigcup_{n \in \nn}\ass \fa_n$ has $t$ elements $\fp_1, \ldots,
\fp_t$. We can therefore relabel the terms in the above-mentioned
minimal primary decompositions so that $\sqrt{\fq_{i,n}} = \fp_i$
for all $i = 1, \ldots, k_n$ and all $n \in \nn$. It is
notationally convenient to define, for any $n \in \nn$ for which
$k_n <t$ and any $j \in \{k_n+1, \ldots, t\}$, an additional
$\fp_j$-primary ideal $\fq_{j,n}$ by $\fq_{j,n} =
f^{-(n'-n)}(\fq_{j,n'})$, where $n'$ is chosen so large that $j
\leq k_{n'}$. We shall then have, for every $n \in \nn$, a (not
necessarily minimal) primary decomposition $\fa_n = \fq_{1,n} \cap
\ldots \cap \fq_{t,n}$ with $\sqrt{\fq_{i,n}}^{[p^n]h} =
\fp_i^{[p^n]h} \subseteq \fq_{i,n}$ for all $i =1, \ldots,t$. (The
extra primary ideals will also have the necessary properties to
ensure that this holds.)

Recall that we have found $c \in R^\circ$ such that  $cr^{p^n}/1
\in (\fa S^{-1}R)^{[p^n]}$ (in the ring of fractions $S^{-1}R$)
for all $n \gg 0$, say for all $n \geq n_0$. Hence, for each $n
\geq n_0$, there exists $s_n \in S$ such that $cs_nr^{p^n} \in
\fa^{[p^n]} \subseteq \fa_n =  \fq_{1,n} \cap \ldots \cap
\fq_{t,n}$. Choose $j \in \{1, \ldots, t\}$. Suppose that there
exists an $m \geq n_0$ such that $cr^{p^m} \not\in \fq_{j,m}$, and
choose the least such $m$. Then $s_m \in \sqrt{\fq_{j,m}} =
\fp_j$: define $s_{\fp_j}$ to be this $s_m$, and note that
$s_{\fp_j}^{hp^n} \in \fq_{j,n}$ for all $n \in \nn$. If, on the
other hand, $cr^{p^n} \in \fq_{j,n}$ for all $n \geq n_0$, set
$s_{\fp_j} = 1$. In both cases, we have $ cs_{\fp_j}^{hp^n}r^{p^n}
\in \fq_{j,n}$ for all $n \geq n_0.$

Set $s :=  s_{\fp_1}\ldots s_{\fp_t} \in S$; we have $
c(rs^h)^{p^n} \in \bigcap_{j=1}^t \fq_{j,n} =\fa_n \subseteq
(\fa^{[p^n]})^*$ for all $n \geq n_0.$ Since $c'$ is a
$p^{m_0}$-weak test element for $R$, we have
$c'(c(rs^h)^{p^n})^{p^{m_0}} \in (\fa^{[p^n]})^{[p^{m_0}]}$ for
all $n \geq n_0$; thus $ c'c^{p^{m_0}}(rs^h)^{p^{n+m_0}} \in
\fa^{[p^{n+m_0}]}$ for all $n \geq n_0.$ Since $c'c^{p^{m_0}} \in
R^{\circ}$, it follows that $rs^h \in \fa^*$, so that $r/1 =
rs^h/s^{h} \in \fa^*S^{-1}R$.
\end{proof}

Our next major aim is to establish that, if a proper ideal $\fA$
of $R^{\infty}$ has a primary decomposition, then the $f$-sequence
$(\fa_n)_{n \in \nn}$ of ideals of $R$ to which it corresponds has
linear growth of primary decompositions, and $\bigcup_{n \in \nn}
\ass \fa_n$ is finite. This will enable us to exploit Theorem
\ref{lg.3}(ii). We need one preparatory lemma.

\begin{lem}
\label{pd.6} Let $\fp \in \Spec (R)$ and let $(\fq_n)_{n \in \nn}$
be an $f$-sequence of $\fp$-primary ideals of $R$. Let $h \in \N$
be such that
 $\fp^{h} \subseteq \fq_0\/.$ Then $\fp^{[p^n]h}=
\sqrt{\fq_{n}}^{[p^n]h} \subseteq \fq_{n}$ for all $n \in \nn$ (so
that, in the language of\/ {\rm \ref{lg.1}}, the $f$-sequence
$(\fq_n)_{n \in \nn}$ has $h$-linear growth of primary
decompositions).
\end{lem}

\begin{proof} By \ref{pd.1a}(i), we have $\fq_0^{[p^n]} \subseteq \fq_{n}
\subseteq \fq_{0}$ for all $n \in \nn$. As $R$ is Noetherian,
there exists $h \in \N$ such that
 $\fp^{h} \subseteq \fq_0\/.$
 Then
 $
  \fq_n \supseteq \fq_0^{[p^n]} \supseteq (\fp^{h})^{[p^n]}
  = \sqrt{\fq_0}^{[p^n]h} = \sqrt{\fq_n}^{[p^n]h}
 $
 for all $n \in \nn $.
\end{proof}

\begin{rmk}
\label{pd.2e} Let $\fp \in \Spec (R)$ and let $(\fq_n)_{n \in
\nn}$ be an $f$-sequence of $\fp$-primary ideals of $R$. Let $n
\in \nn$. By \ref{pd.1a}(i), we have $(\fq_0^{[p^n]})^F \subseteq
\fq_n$. Let $\phantom{\fa}^e$ and $\phantom{\fa}^c$ stand for
extension and
 contraction with respect to the natural ring homomorphism
$R \lra R_{\fp}.$ Then $((\fq_0^{[p^n]})^F)^{ec} \subseteq \fq_n$,
since $\fq_n$ is $\fp$-primary. A simple example quickly shows
that we cannot expect equality here: let $K$ be a field of
characteristic $p$, take $R = K[X]$, the polynomial ring in one
indeterminate, and take $\fq_n = (X)$ for all $n \in \nn$.
\end{rmk}

\begin{thm}
\label{pd.7} If the proper ideal $\fA$ of $R^\infty$ has a primary
decomposition, then the $f$-sequence $(\fa_n)_{n \in \nn}$ of
ideals of $R$ corresponding to $\fA$ has linear growth of primary
decompositions, and\/ $\bigcup_{n \in \nn}\ass \fa_n$ is a finite
set.
\end{thm}

\begin{proof}  By Theorem \ref{pd.5}, there are prime ideals
$\fp_1,\ldots,\fp_t$ of $R$ such that the following is true: each
$\fa_n$ has a primary decomposition
    $\fa_n = \fq_{1,n} \cap \dots \cap \fq_{t,n}$ such that
    $\fq_{k,n}$ is $\fp_k$-primary
       for all $k = 1,
\ldots, t$ and all $n \in \nn$, and
      $(\fq_{k,n})_{n \in \nn}$ is an
       $f$-sequence for all $k = 1,
\ldots, t\/.$

Now there exist positive integers $h_1, \ldots, h_t$ such that
$\sqrt{\fq_{k,0}}^{h_k} \subseteq \fq_{k,0}$ for all $k = 1,
\ldots, t$.
 By Lemma \ref{pd.6}, we have
 $
  \sqrt{\fq_{k,n}}^{[p^n]h_k} \subseteq \fq_{k,n}$
  for all $n \in \nn$ and all $k = 1,
\ldots, t.$
 Hence, with $h := \max \{h_1, \ldots, h_t\}$, we have
$
  \sqrt{\fq_{k,n}}^{[p^n]h} \subseteq \fq_{k,n}$
  for all $n \in \nn$ and all $k = 1,
\ldots, t.$  The final claim follows also from Theorem \ref{pd.5}.
\end{proof}

In view of Theorems \ref{lg.3}(ii) and
\ref{pd.7}, we are very interested in
finding primary decompositions of proper ideals of $R^{\infty}$.
However, in Theorem \ref{pd.4} we showed that $R^{\infty}$ is only
Noetherian in rather uninteresting cases, and so the existence of
a primary decomposition for a proper ideal of $R^{\infty}$ would
be a bonus that we should not expect in all cases. Indeed,
Example \ref{ex.8}
below is of an $f$-sequence $(\a_n)_{n \in \nn}$ of proper ideals
(in a $2$-dimensional regular ring $R$)
which has linear growth of primary decompositions, but for which
the set $\bigcup_{n \in \nn}\ass \fa_n$ is infinite (so that, by Theorem \ref{pd.7},
the associated ideal of $R^{\infty}$ cannot have a primary
decomposition).

\begin{ex}
\label{ex.8} Let $K$ be an infinite field of prime characteristic
$p$, let $(\lambda_j)_{j \in \N}$ be a sequence of distinct
elements of $K$, and let $(t_j)_{j \in \N}$ be a sequence of
positive integers. Let $R = K[X,Y]$, where $X$ and $Y$ are
independent indeterminates. Let $l$ be an integer such that $2
\leq l \leq p$.

Let $(\fq_{0,n})_{n \in \nn}$ be the $f$-sequence of ideals of $R$
for which $\fq_{0,n} = (X)$ for all $n \in \nn$: see
\ref{pd.2}(i).

Note that $R$ is regular, so that, if $\fa$ is an arbitrary ideal
of $R$, then $(\fa^{[p^n]})_{n \in \nn}$ is the $f$-sequence
$((\fa^{[p^n]})^F)_{n \in \nn}$, since every ideal of $R$ is
tightly closed (see \ref{pd.0}(iii)). With this and
\ref{pd.1a}(iii) in mind, we let, for each $j \in \N$,
$(\fq_{j,n})_{n \in \nn}$ be the $f$-sequence of
$(X,(Y-\lambda_j))$-primary ideals of $R$ given by
$$
\fq_{j,n} = \begin{cases}
             \big(X^{lp^{n-j}},(Y-\lambda_j)^{t_jp^n}\big) =
             \big(X^{l},(Y-\lambda_j)^{t_jp^j}\big)^{[p^{n-j}]}& \text{if $n \geq j$,}   \\
         f^{-(j-n)}(\fq_{j, j}) & \text{if $n < j$.}
\end{cases}
$$
Let $(\fa_{n})_{n \in \nn}$ be the $f$-sequence of ideals of $R$
given by $\fa_n := \bigcap_{j\in \nn} \fq_{j,n}$ for all $n \in
\nn$; see \ref{pd.2}(v). Note that, for $m,j \in \nn$ with $j >
m$, we have $ f^{j-m}(X) = X^{p^{j-m}} = X^lX^{p^{j-m}-l} \in
(X^{l},(Y-\lambda_j)^{t_jp^j}) = \fq_{j,j}$; hence $\fq_{0,m} =
(X) \subseteq f^{-(j-m)}(\fq_{j,j}) = \fq_{j,m}$ for all $j > m$.
Thus
$$
\fa_{m} = (X) \cap \big(X^{lp^{m-1}},(Y-\lambda_1)^{t_1p^m}\big)
\cap \big(X^{lp^{m-2}},(Y-\lambda_2)^{t_2p^m}\big) \cap \ldots
\cap \big(X^{l},(Y-\lambda_m)^{t_mp^m}\big).
$$
It therefore follows that $ \fa_m = \fq_{0,m} \cap \fq_{1,m} \cap
\ldots \cap \fq_{m,m} $ is a primary decomposition (in which the
radicals of the primary terms are all different), for each $m \in
\nn$. We show next, by induction on $m$, that these primary
decompositions are all minimal; it is clear that $\fa_0 =
\fq_{0,0}$ is a minimal primary decomposition, and so we suppose
now that $m > 0$ and make the inductive hypothesis that $
\fa_{m-1} = \fq_{0,m-1} \cap \fq_{1,m-1} \cap \ldots \cap
\fq_{m-1,m-1} $ is a minimal primary decomposition.

Since $(\fa_{n})_{n \in \nn}$ is an $f$-sequence, $\ass \a_{m-1}
\subseteq \ass \a_{m}$ by \ref{pd.1a}(iv). Therefore, none of
$\fq_{0,m}, \ldots,\fq_{m-1,m}$ can be omitted from the primary
decomposition $\fa_m = \bigcap_{k=1}^m\fq_{k,m}$. We can then
conclude that this primary decomposition is minimal simply by
observing that
$$
X(Y-\lambda_1)^{t_1p^m}\ldots(Y-\lambda_{m-1})^{t_{m-1}p^m} \in
\left(\fq_{0,m} \cap \fq_{1,m} \cap \ldots \cap \fq_{m-1,m}\right)
\setminus \fq_{m,m}.
$$

A consequence of this inductive argument is that $(\fa_{n})_{n \in
\nn}$ is an $f$-sequence of ideals of $R$ for which
$$
\ass \a_n = \left\{ (X), (X, (Y-\lambda_1) ), \ldots, (X,
(Y-\lambda_n) ) \right\} \quad \mbox{~for all~} n \in \nn.
$$
(Of course, in the case when $n = 0$, this statement is to be
interpreted as `$\ass \a_0 = \left\{ (X) \right\}$'.)
It therefore follows from Theorem \ref{pd.5} that
the ideal of $R^{\infty}$ to which $(\fa_{n})_{n \in \nn}$
corresponds does not have a primary decomposition. However, if the
sequence $(t_j)_{j \in \N}$ is chosen so that it is bounded, by an
integer $t$ say, then it is straightforward to check that
$\sqrt{\fq_{j,n}}^{[p^n]t} \subseteq \fq_{j,n}$ for all $n,j \in
\nn$ with $j \leq n$; thus $(\fa_{n})_{n \in \nn}$ does have
linear growth of primary decompositions.
\end{ex}

We revert now to the general situation of our standard hypotheses.
Our next result will show that the primary components of the terms
$\fa_n$ in an $f$-sequence $(\fa_n)_{n \in\nn}$ corresponding to
the minimal prime ideals of the $f$-sequence (see \ref{pd.1a}(v))
never present any obstacle to the $f$-sequence's having linear
growth of primary decompositions.

\begin{prop}
\label{pd.2b} Let $(\fa_n)_{n \in\nn}$ be an $f$-sequence of
proper ideals of $R$, and let $\fp_1, \ldots, \fp_l$ be the
minimal primes of this $f$-sequence (see\/ {\rm \ref{pd.1a}(v)}).
For each $n \in \nn$ and each $k \in \{1, \ldots,l\}$, let
$\fq_{k,n}$ be the (uniquely determined) $\fp_k$-primary component
of $\fa_n$. Let $k_1, \ldots,k_t$ be integers with $1 \leq k_1 <
k_2 < \ldots < k_t \leq l$. Then
$
\left(\fq_{k_1,n} \cap \fq_{k_2,n} \cap \ldots \cap
\fq_{k_t,n}\right)_{n
 \in \nn}
$
is an $f$-sequence of ideals of $R$ which has linear growth of
primary decompositions, and the ideal of $R^{\infty}$ to which it
corresponds has a primary decomposition.
\end{prop}

\begin{proof} Set $S := R \setminus \bigcup_{j=1}^t \fp_{k_j}$, a
multiplicatively closed subset of $R$. Let $\phantom{\fa}^e$ and
$\phantom{\fa}^c$ stand for extension and
 contraction with respect to the natural ring homomorphism
$R \lra S^{-1}R\/.$
 By \ref{pd.2a}, the sequence
$(\fa_n^{ec})_{n \in\nn}$ is an $f$-sequence of ideals of $R\/.$
But $ \fa_n^{ec} = \fq_{k_1,n} \cap \fq_{k_2,n} \cap \ldots \cap
\fq_{k_t,n}$ for all $n \in \nn.$ As particular cases, we see
that, for each $j = 1, \ldots, t$, the sequence
$(\fq_{k_j,n})_{n\in \nn}$ is an $f$-sequence. It therefore
follows from Theorem \ref{pd.5} that the ideal of $R^{\infty}$ to
which the $f$-sequence $(\fa_n^{ec})_{n \in\nn}$ corresponds has a
primary decomposition. Hence $(\fa_n^{ec})_{n \in\nn}$ has linear
growth of primary decompositions, by Theorem \ref{pd.7}.
\end{proof}

\begin{cor}
\label{lg.11} Let $(\fa_n)_{n \in\nn}$ be an $f$-sequence of
proper ideals of $R$ with the property that each $\fa_n$ $(n \in
\nn)$ has no embedded prime. Then $(\fa_n)_{n \in\nn}$ has linear
growth of primary decompositions, and the corresponding ideal of
$R^{\infty}$ has a primary decomposition.
\end{cor}

Thus the problems, in showing that a given $f$-sequence
$(\fa_n)_{n \in\nn}$ of proper ideals of $R$ has linear growth of
primary decompositions, rest entirely with the embedded primary
components of the $\fa_n$. None of these is ever uniquely
determined, and the issues revolve around whether or not it is
possible to make appropriate choices for these embedded primary
components. In (the proof of) Theorem \ref{pd.7} above, we saw
that, if the ideal of $R^{\infty}$ corresponding to the
$f$-sequence $(\fa_n)_{n \in\nn}$ has a primary decomposition,
then there are natural choices for the above-mentioned embedded
primary components that satisfy the conditions necessary for
linear growth of primary decompositions.

In the next section, we shall present some examples of proper
ideals in non-Noetherian perfect closures that do have primary
decompositions.

The results of this section provide a strategy for attempting to
show that, for a given proper ideal $\fa$ of $R$, tight closure
commutes with localization with respect to an arbitrary
multiplicatively closed subset of $R$: if we can find an $f$-sequence
$(\fa_n)_{n \in\nn}$ of ideals of $R$ such that $\fa_n$ is trapped over
$\fa^{[p^n]}$ for all $n \in \nn$, if $R$ has a $p^{m_0}$-weak test element
(for some $m_0 \in \nn$) and if the ideal of $R^{\infty}$ to
which $(\fa_n)_{n \in\nn}$ corresponds has a primary decomposition, then
it follows from Theorem \ref{pd.7} that $(\fa_n)_{n \in\nn}$ has linear growth of
primary decompositions and that $\bigcup_{n \in \nn}\ass \fa_n$ is a finite
set, and it then follows from Theorem \ref{lg.3}(ii) that
$\fa^*S^{-1}R =
(\fa S^{-1}R)^*$ for every multiplicatively closed subset $S$ of
$R$. We should perhaps mention at this point that, by-and-large,
we have only managed to get this strategy to succeed in situations where it
had already been proved that tight closure commutes with localization. However, our
next result shows that the
above-mentioned hypotheses needed for the strategy to work actually
ensure that the $f$-sequence $\big((\fa^{[p^n]})^*\big)_{n \in \nn}$ of tight
closures of the Frobenius powers of $\fa$ (see Lemma \ref{pd.1c}) also has
linear growth of primary decompositions, and that $\bigcup_{n \in
\nn}\ass (\fa^{[p^n]})^*$ is finite.

This is interesting because
M. Katzman's approach in \cite{Katzm96} to the localization problem
for tight closure led to the following
question: is it the case that, for every ideal $\fb$ of $R$, the set
$\bigcup_{n \in \nn}\ass (\fb^{[p^n]})^*$
has only finitely many maximal elements? Our strategy outlined above, used
in conjunction with Theorem \ref{lg.33} below, will enable us to conclude in
\S\ref{ex} that $\bigcup_{n \in \nn}\ass (\fb^{[p^n]})^*$ is actually
a finite set in
several cases where it is known that, for $\fb$, tight closure commutes with
localization with respect to an arbitrary multiplicatively closed subset of $R$.

\begin{thm}
\label{lg.33} Let $\fa$ be a proper ideal of $R$.  Suppose that
$(\fa_n)_{n \in\nn}$ is an $f$-sequence of ideals of $R$ which is
such that $\fa^{[p^n]} \subseteq \fa_n \subseteq(\fa^{[p^n]})^*$
for all $n \in \nn$, and that the ideal $\fA$ of $R^{\infty}$ to
which $(\fa_n)_{n \in\nn}$ corresponds has a primary
decomposition, so that, by\/ {\rm \ref{pd.7}}, the set\/
$\bigcup_{n \in \nn}\ass \fa_n$ is finite. Let the members of the
latter set be $\fp_1, \ldots, \fp_t$, and let the corresponding
prime ideals of $R^{\infty}$ (see\/ {\rm \ref{pd.2}(i)}) be
$\fP_1, \ldots, \fP_t$ respectively.

For each $i = 1, \ldots, t$, let $\phantom{\fa}^{e_i}$ and
$\phantom{\fa}^{c_i}$ stand for extension and
 contraction with respect to the natural ring homomorphism
$R \lra R_{\fp_i}$.

\begin{enumerate}
\item There exists $h \in \N$ such that, for each $n \in \nn$,
$$
\fa_n = \bigcap_{i=1}^t \big((\fa_n +
\fp_i^{h[p^n]})^F\big)^{e_ic_i}, \quad \mbox{~where~} \big((\fa_n
+ \fp_i^{h[p^n]})^F\big)^{e_ic_i} \mbox{~is~} \fp_i\mbox{-primary
for~} i = 1, \ldots,t,
$$
is a primary decomposition in which each primary component is $F$-closed.
\item
Assume that $R$ has a $p^{m_0}$-weak test element, for
some $m_0 \in \nn$. Then, with $h$ as in\/ {\rm (i)},

\begin{enumerate}
\item for each $n \in \nn$,
$$
(\fa^{[p^n]})^* = \bigcap_{i=1}^t \big(((\fa +
\fp_i^h)^{[p^n]})^*\big)^{e_ic_i}, \quad \mbox{~where~} \big(((\fa
+ \fp_i^h)^{[p^n]})^*\big)^{e_ic_i} \mbox{~is~}
\fp_i\mbox{-primary for~} i = 1, \ldots,t,
$$
is a primary decomposition in which each primary component is $F$-closed,
and the ideal of $R^{\infty}$ to which the $f$-sequence
$\big((\fa^{[p^n]})^*\big)_{n \in \nn}$ corresponds has a primary decomposition;

\item the $f$-sequence $\big((\fa^{[p^n]})^*\big)_{n \in \nn}$
has $h$-linear growth of primary decompositions, and the set\/
$\bigcup_{n \in \nn} \ass (\fa^{[p^n]})^*$ is finite; and
\item $\bigcup_{n \in \nn} \ass (\fa^{[p^n]})^* \subseteq
\bigcup_{n \in \nn} \ass \fa_n$.
\end{enumerate}
\item
If $R$ has a locally stable $p^{m_0}$-weak test element, for
some $m_0 \in \nn$, then, for each $n \in \nn$, and with $h$ as in\/ {\rm (i)},
$$
(\fa^{[p^n]})^* = \bigcap_{i=1}^t \big((((\fa +
\fp_i^h)^{[p^n]})^{e_i})^*\big)^{c_i}, \quad
\mbox{~where~}\big((((\fa + \fp_i^h)^{[p^n]})^{e_i})^*\big)^{c_i}
\mbox{~is~} \fp_i\mbox{-primary for~} i = 1, \ldots,t,
$$
is a primary decomposition in which each primary component is tightly closed.
\end{enumerate}
\end{thm}

\begin{proof} By Theorems \ref{pd.5} and \ref{pd.7}, each
$\fa_n$ has a primary decomposition
    $\fa_n = \fq_{1,n} \cap \dots \cap \fq_{t,n}$ such that
          $(\fq_{i,n})_{n \in \nn}$ is an
       $f$-sequence of $\fp_i$-primary ideals of $R$ for all $i = 1,
\ldots, t$, and there exists $h \in \N$ such that $
  \sqrt{\fq_{i,n}}^{h[p^n]} = \fp_i^{h[p^n]} \subseteq \fq_{i,n} = \fq_{i,n}^F$
  for all $n \in \nn$ and all $i = 1,
\ldots, t.$

(i) Now $\fa_n \subseteq \fa_n + \fp_i^{h[p^n]}
 \subseteq \fq_{i,n}$,
so that $$\fa_n
\subseteq \big((\fa_n
+ \fp_i^{h[p^n]})^F\big)^{e_ic_i} \subseteq (\fq_{i,n}^F)^{e_ic_i} =
\fq_{i,n}^{e_ic_i} = \fq_{i,n} \quad \mbox{~for all~} i =
1, \ldots, t.$$ Hence $\fa_n = \bigcap_{i = 1}^t
\big((\fa_n
+ \fp_i^{h[p^n]})^F\big)^{e_ic_i}$; moreover $$\fp_i^{h[p^n]}  \subseteq
\big(\fa_n
+ \fp_i^{h[p^n]}\big)^{e_ic_i}  \subseteq
\big((\fa_n
+ \fp_i^{h[p^n]})^F\big)^{e_ic_i} \subseteq \fp_i \quad
\mbox{~for all~} i = 1, \ldots, t,$$ so that $\big((\fa_n
+ \fp_i^{h[p^n]})^F\big)^{e_ic_i}$ is $\fp_i$-primary. It follows
that
$$
\fa_n = \bigcap_{i=1}^t \big((\fa_n +
\fp_i^{h[p^n]})^F\big)^{e_ic_i}, \quad \mbox{~where~} \big((\fa_n
+ \fp_i^{h[p^n]})^F\big)^{e_ic_i} \mbox{~is~} \fp_i\mbox{-primary
for~} i = 1, \ldots,t,
$$
is a primary decomposition, and parts (iii) and (i) of Lemma \ref{bk.3} show that
each primary component in this decomposition is $F$-closed.

(ii) Let $c$ be a $p^{m_0}$-weak test element for $R$.

It is clear that $(\fa^{[p^n]})^* \subseteq
\bigcap_{i=1}^t ((\fa + \fp_i^h)^{[p^n]})^* \subseteq
\bigcap_{i=1}^t \big(((\fa + \fp_i^h)^{[p^n]})^* \big)^{e_ic_i}$.
Note that, for each $i \in \{1, \ldots,t\}$, we have $
\fp_i^{[p^n]h} \subseteq (\fa + \fp_i^h)^{[p^n]}  \subseteq
\fp_i$, so that
$$
\fp_i^{[p^n]h} \subseteq ((\fa + \fp_i^h)^{[p^n]})^* \subseteq
\big(((\fa + \fp_i^h)^{[p^n]})^* \big)^{e_ic_i} \subseteq
\big(\fp_i^*\big)^{e_ic_i} = \fp_i^{e_ic_i} = \fp_i.
$$
It therefore follows that $\big(((\fa + \fp_i^h)^{[p^n]})^*
\big)^{e_ic_i}$ is $\fp_i$-primary. We show next that
$\bigcap_{i=1}^t
\big(((\fa + \fp_i^h)^{[p^n]})^* \big)^{e_ic_i} \subseteq (\fa^{[p^n]})^*$.

So let $r \in \bigcap_{i=1}^t \big(((\fa + \fp_i^h)^{[p^n]})^*
\big)^{e_ic_i}.$ Then, for each $i \in \{1, \ldots, t\}$, there
exists $s_i \in R \setminus \fp_i$ such that $s_ir \in ((\fa +
\fp_i^h)^{[p^n]})^*$; therefore, for all $m \geq m_0$, we have
$c(s_ir)^{p^m} \in ((\fa + \fp_i^h)^{[p^n]})^{[p^m]}$. Therefore
$$
cs_i^{p^m}r^{p^m} \in \fa^{[p^{n+m}]} + \fp_i^{h[p^{n+m}]} \subseteq
\fa_{n+m}  + \fp_i^{h[p^{n+m}]}  \subseteq
\big(\fa_{n+m}  + \fp_i^{h[p^{n+m}]}\big)^F \quad \mbox{~for all~} m \geq m_0.
$$
Hence, for all $i = 1, \ldots, t$, we have $cr^{p^m} \in
\big((\fa_{n+m}  + \fp_i^{h[p^{n+m}]})^F\big)^{e_ic_i}$
for all $m \geq m_0$. One can now
use part (i) to deduce that $cr^{p^m} \in \fa_{n+m} \subseteq
(\fa^{[p^{n+m}]})^*$ for all $m \geq m_0$.

We use again the fact that $c$ is a $p^{m_0}$-weak test element to deduce that
$$
c(cr^{p^m})^{p^k} \in (\fa^{[p^{n+m}]})^{[p^k]} = \fa^{[p^{n+m+k}]} \quad
\mbox{~for all~} m,k \geq m_0.
$$
Take $k = m_0$ to see that $c^{1 + p^{m_0}}r^{p^{m+m_0}} \in
(\fa^{[p^{n}]})^{[p^{m+m_0}]}$ for all $m \geq m_0$. Since $c^{1 + p^{m_0}} \in
R^{\circ}$, it follows that $r \in (\fa^{[p^{n}]})^*$. Therefore
$$
(\fa^{[p^n]})^* = \bigcap_{i=1}^t \big(((\fa +
\fp_i^h)^{[p^n]})^*\big)^{e_ic_i}, \quad \mbox{~where~} \big(((\fa
+ \fp_i^h)^{[p^n]})^*\big)^{e_ic_i} \mbox{~is~}
\fp_i\mbox{-primary for~} i = 1, \ldots,t,
$$
is a primary decomposition. Since a tightly-closed ideal is $F$-closed, it again
follows from parts (iii) and (i) of Lemma \ref{bk.3} that
each primary component in this decomposition is $F$-closed.

By Lemmas \ref{pd.1c} and \ref{pd.2a}, for each $i = 1, \ldots, t$,
the sequence $\left(\big(((\fa
+ \fp_i^h)^{[p^n]})^*\big)^{e_ic_i}\right)_{n\in\nn}$ is an $f$-sequence
of $\fp_i$-primary ideals. It
follows from Theorem \ref{pd.5} that the ideal of $R^{\infty}$
to which the $f$-sequence $\left((\fa^{[p^n]})^*\right)_{n\in\nn}$ corresponds
has a primary decomposition; all the remaining claims in part (ii) are now clear.

(iii)  Let $c$ be a locally stable $p^{m_0}$-weak test element for $R$.

Arguments similar to those used in the above proof of part (ii) will show that
$$(\fa^{[p^n]})^* \subseteq \bigcap_{i=1}^t \big((((\fa +
\fp_i^h)^{[p^n]})^{e_i})^*\big)^{c_i}$$ and that $\big((((\fa +
\fp_i^h)^{[p^n]})^{e_i})^*\big)^{c_i}$ is $\fp_i$-primary for each
$i = 1, \ldots, t$.

Now let $r \in \bigcap_{i=1}^t \big((((\fa +
\fp_i^h)^{[p^n]})^{e_i})^*\big)^{c_i}$. Let $j \in \{1, \ldots, t\}$. Then
$r/1 \in (((\fa +
\fp_j^h)^{[p^n]})^{e_j})^*$; since $c/1$ is a $p^{m_0}$-weak test element
for $R_{\fp_j}$, we see that, for all $m \geq m_0$,
\begin{align*}
cr^{p^m}/1 \in (((\fa +
\fp_j^h)^{[p^n]})^{e_j})^{[p^m]} &= ((\fa +
\fp_j^h)^{[p^{m+n}]})^{e_j}
= (\fa^{[p^{m+n}]} +
\fp_j^{h[p^{m+n}]})^{e_j} \\
&\subseteq  (\fa_{m+n} +
\fp_j^{h[p^{m+n}]})^{e_j} \subseteq ((\fa_{m+n} +
\fp_j^{h[p^{m+n}]})^F)^{e_j}.
\end{align*}
Therefore $cr^{p^m} \in
\big((\fa_{n+m}  + \fp_j^{h[p^{n+m}]})^F\big)^{e_jc_j}$
for all $m \geq m_0$. One can now
use part (i) once again to deduce that $cr^{p^m} \in \fa_{n+m} \subseteq
(\fa^{[p^{n+m}]})^*$ for all $m \geq m_0$, and proceed as in the proof of part (ii)
to conclude that $r \in (\fa^{[p^n]})^*$.

Thus
$$
(\fa^{[p^n]})^* = \bigcap_{i=1}^t \big((((\fa +
\fp_i^h)^{[p^n]})^{e_i})^*\big)^{c_i}, \quad
\mbox{~where~}\big((((\fa + \fp_i^h)^{[p^n]})^{e_i})^*\big)^{c_i}
\mbox{~is~} \fp_i\mbox{-primary for~} i = 1, \ldots,t,
$$
is a primary decomposition, and this time one can use Lemma \ref{bk.3}(ii) to see
that $\big((((\fa +
\fp_i^h)^{[p^n]})^{e_i})^*\big)^{c_i}$ is tightly closed (for each
$i = 1, \ldots, t$).
\end{proof}

The next theorem draws together results that now follow, in the presence
of additional hypotheses, from the original hypothesis of Smith
and Swanson that the Frobenius powers of the proper ideal
$\fa$ of $R$ have linear growth of primary decompositions.

\begin{thm}
\label{lg.4} Suppose that the Frobenius powers of the proper ideal
$\fa$ of $R$ have linear growth of primary decompositions, and
that\/ $\bigcup_{n \in \nn} \ass \fa^{[p^n]}$ is finite. Then
\begin{enumerate}
\item the ideal $\fa R^{\infty}$ of $R^{\infty}$ to which the $f$-sequence
$\big((\fa^{[p^n]})^F\big)_{n \in \nn}$ corresponds has a primary decomposition;
\item the $f$-sequence $\big((\fa^{[p^n]})^F\big)_{n \in \nn}$
has linear growth of primary decompositions;
\item $\bigcup_{n \in \nn} \ass (\fa^{[p^n]})^F \subseteq
\bigcup_{n \in \nn} \ass  \fa^{[p^n]}$, and so the set\/ $\bigcup_{n \in \nn}
\ass (\fa^{[p^n]})^F$ is finite;
\end{enumerate}
furthermore, if $R$ has a $p^{m_0}$-weak
test element, for some $m_0 \in \nn$, then, in addition,
\begin{enumerate}
\setcounter{enumi}{3}
\item $\fa^*S^{-1}R = (\fa S^{-1}R)^*$ for every multiplicatively closed subset
$S$ of $R$;
\item the ideal of $R^{\infty}$ to which the $f$-sequence
$\big((\fa^{[p^n]})^*\big)_{n \in \nn}$ corresponds has a primary
decomposition, so that the $f$-sequence
$\big((\fa^{[p^n]})^*\big)_{n \in \nn}$ has linear growth of
primary decompositions; and
\item $\bigcup_{n \in \nn} \ass (\fa^{[p^n]})^* \subseteq
\bigcup_{n \in \nn} \ass  (\fa^{[p^n]})^F$, so that the set\/
$\bigcup_{n \in \nn} \ass (\fa^{[p^n]})^*$ is finite.
\end{enumerate}
\end{thm}

\begin{proof}
The claims in parts (i), (ii) and (iii) follow from Theorem
\ref{lg.2}. Since $\fa^{[p^n]} \subseteq (\fa^{[p^n]})^F \subseteq
(\fa^{[p^n]})^*$ for all $n \in \nn$, all the other claims
therefore follow from Theorems \ref{lg.3}(ii) and \ref{lg.33}(ii).
\end{proof}

We end this section by drawing together in one theorem the conclusions
of this section which together provide a strategy for attempting to show,
for a given proper ideal $\fa$ of $R$, that tight closure commutes with localization
at an arbitrary multiplicatively closed subset of $R$ and that
$\bigcup_{n \in \nn}
\ass (\fa^{[p^n]})^*$ is finite.

\begin{thm}
\label{lg.3a} Let $\fa$ be a proper ideal of $R$. Suppose that
$(\fa_n)_{n \in\nn}$ is an $f$-sequence of ideals of $R$ which is
such that $\fa^{[p^n]} \subseteq \fa_n \subseteq(\fa^{[p^n]})^*$
for all $n \in \nn$. Suppose that the ideal $\fA$ of $R^{\infty}$
to which $(\fa_n)_{n \in\nn}$ corresponds has a primary
decomposition. Then

\begin{enumerate}
\item the $f$-sequence $(\fa_n)_{n \in\nn}$
has linear growth of primary decompositions;
\item the set\/ $\bigcup_{n \in
\nn}\ass \fa_n$ is finite;
\end{enumerate}
furthermore, if $R$ has a $p^{m_0}$-weak
test element, for some $m_0 \in \nn$, then, in addition,
\begin{enumerate}
\setcounter{enumi}{2}
\item $\fa^*S^{-1}R = (\fa S^{-1}R)^*$ for every multiplicatively closed
subset $S$ of $R$;
\item the ideal of $R^{\infty}$ to which the $f$-sequence
$\big((\fa^{[p^n]})^*\big)_{n \in \nn}$ corresponds has a primary
decomposition, so that the $f$-sequence
$\big((\fa^{[p^n]})^*\big)_{n \in \nn}$ has linear growth of
primary decompositions; and
\item $\bigcup_{n \in \nn} \ass (\fa^{[p^n]})^* \subseteq
\bigcup_{n \in \nn} \ass  \fa_n$, so that the set\/ $\bigcup_{n
\in \nn} \ass (\fa^{[p^n]})^*$ is finite.
\end{enumerate}
\end{thm}

\begin{proof} The claims in (i) and (ii) follow from Theorem \ref{pd.7}. The
claim in (iii) then follows from Theorem \ref{lg.3}(ii), while
those in (iv) and (v) follow from Theorem \ref{lg.33}(ii).
\end{proof}

\section{\sc Applications of the strategy}
\label{ex}

Throughout this section, $R$ will denote a commutative Noetherian
ring of prime characteristic $p$. The aim of this section is to
give several examples of choices of ($R$ and) the proper ideal
$\fa$ of $R$ for which the set $\bigcup_{n \in \nn} \ass
(\fa^{[p^n]})^*$ is a finite set. We use either the strategy of
Theorem \ref{lg.4} or that of Theorem \ref{lg.3a}, but in most of
the examples presented below, it is already known that
$\fa^*S^{-1}R = (\fa S^{-1}R)^*$ for every multiplicatively closed
subset $S$ of $R$. Our first application uses Theorem \ref{lg.4}.

\begin{prop}
\label{me.4} Let $\fa$ be a proper ideal of $R$ such that, for every $n \in \nn$,
all the associated primes of $\fa^{[p^n]}$ are minimal. Then
\begin{enumerate}
\item the Frobenius powers of
$\fa$ have linear growth of primary decompositions;
\item the ideal $\fa R^{\infty}$ of $R^{\infty}$ has a primary decomposition;
\item the $f$-sequence $\big((\fa^{[p^n]})^F\big)_{n \in \nn}$
has linear growth of primary decompositions;
\item $\bigcup_{n \in \nn} \ass (\fa^{[p^n]})^F$ is equal to the set of minimal
primes of $\fa$;
\end{enumerate}
furthermore, if $R$ has a $p^{m_0}$-weak
test element, for some $m_0 \in \nn$, then, in addition,
\begin{enumerate}
\setcounter{enumi}{4}
\item $\fa^*S^{-1}R = (\fa S^{-1}R)^*$ for every multiplicatively closed subset
$S$ of $R$;
\item the ideal of $R^{\infty}$ to which the $f$-sequence
$\big((\fa^{[p^n]})^*\big)_{n \in \nn}$ corresponds has a primary
decomposition, and so the $f$-sequence
$\big((\fa^{[p^n]})^*\big)_{n \in \nn}$ has linear growth of
primary decompositions; and
\item the set\/
$\bigcup_{n \in \nn} \ass (\fa^{[p^n]})^*$ is (finite and)
equal to the set of minimal
primes of $\fa$.
\end{enumerate}
\end{prop}

\begin{proof} Let $\fp_1, \ldots,
\fp_t$ be the minimal prime ideals of $\fa$.
For each $i = 1, \ldots, t$, let $\phantom{\fa}^{e_i}$ and
$\phantom{\fa}^{c_i}$ stand for extension and
 contraction with respect to the natural ring homomorphism
$R \lra R_{\fp_i}$.

Now $\fa^{[p^n]} = \bigcap_{i=1}^t (\fa^{[p^n]})^{e_ic_i}$ is the
unique minimal primary decomposition of $\fa^{[p^n]}$ (for each $n \in \nn$).
Choose $h \in \N$ such that $\fp_i^h \subseteq \fa^{e_ic_i}$ for all
$i = 1, \ldots, t$. Let $i \in \{1, \ldots, t\}$ and $r \in \fp_i^h$. Then
there exists $s_i \in R \setminus \fp_i$ such that $s_ir \in \fa$, so that
$s_i^{p^n}r^{p^n} \in \fa^{[p^n]}$ and $r^{p^n} \in (\fa^{[p^n]})^{e_ic_i}$
(for all $n \in \nn$). This shows that the Frobenius powers of
$\fa$ have $h$-linear growth of primary decompositions.

The remaining claims follow from Theorem \ref{lg.4} (and the fact
that a prime ideal of $R$ is tightly closed).
\end{proof}

\begin{ex}
\label{me.5} If $R$ is a Cohen-Macaulay ring and $\fa$ is an ideal
generated by a regular sequence $(r_i)_{i=1}^k$, then, for all $n
\in \nn$, the $n$-th Frobenius power $\fa^{[p^n]}$ is generated by
the regular sequence $(r_i^{p^n})_{i=1}^k$, and so is unmixed;
thus the initial hypothesis of \ref{me.4} is satisfied, and
conclusions (i), (ii), (iii), (iv) of \ref{me.4} hold. If, in
addition, $R$ has a $p^{m_0}$-weak test element, for some $m_0 \in
\nn$, then the other three conclusions of \ref{me.4} also hold.
\end{ex}

Some of the results in this section will concern situations where we can, in some
sense, `approximate' $R$ by a regular commutative Noetherian ring
of characteristic $p$.

\begin{rmk}
\label{ex.2} Suppose that $R$ is regular.
\begin{enumerate}
\item As the
Frobenius homomorphism $f : R \lra R$ is flat (see \cite{Kunz69}),
every ideal of $R$ is $F$-closed (and tightly closed).
\item It also follows from the fact that $f$ is flat that,
if $\fq$ is a $\fp$-primary ideal of $R$, then so too is
$\fq^{[p^n]}$ for all $n \in \nn$, and that if $\fa = \fq_1 \cap
\ldots \cap \fq_t$ is a minimal primary decomposition of the
proper ideal $\fa$ of $R$, then $\fa^{[p^n]} = \fq_1^{[p^n]} \cap
\ldots \cap \fq_t^{[p^n]}$ is also a minimal primary decomposition
of $\fa^{[p^n]}$, for all $n \in \nn$.
\item It is even the case that, if $\fq$ is a $\fp$-primary
irreducible ideal of $R$, then $\fq^{[p^n]}$ is also irreducible,
for each $n \in \nn$. To see this, let $F_R$ denote the Frobenius
functor $R'\otimes_R(\: {\scriptscriptstyle \bullet} \:)$ on the
category of $R$-modules and $R$-homomorphisms, where $R'$ is $R$
considered as a left $R$-module in the usual manner and as a right
$R$-module via the Frobenius homomorphism.

Since $E_R(R/\fq) \cong E_R(R/\fp)$, it follows that $R/\fq^{[p]}
\cong F(R/\fq)$ is isomorphic to a submodule of $F(E_R(R/\fp))$;
however, $F(E_R(R/\fp)) \cong E_R(R/\fp)$, by \cite[Proposition
1.5]{58}, and so $E_R\big(R/\fq^{[p]}\big)$ must be isomorphic to
an injective submodule of the indecomposable injective $R$-module
$E_R(R/\fp)$; therefore $E_R\big(R/\fq^{[p]}\big) \cong
E_R(R/\fp)$ and $\fq^{[p]}$ is irreducible.
\end{enumerate}
\end{rmk}

\begin{prop}
\label{ex.3} Suppose that $R$ is regular. Then each finitely
generated proper ideal of $R^{\infty}$ has a primary
decomposition.
\end{prop}

\begin{proof}
Let $\fA$ be a finitely generated proper ideal in $R^\infty\/.$
Let $(\fa_n)_{n \in\nn}$
 be the $f$-sequence corresponding to $\fA\/.$
As every ideal in $R$ is $F$-closed (see \ref{ex.2}(i)),
 it follows from Theorem \ref{ex.1} that there exists
$k \in \nn$ such that $\fa_{k+n} = \fa_k^{[p^n]}$ for all $n \in
\nn$.

 Consider a minimal primary decomposition
 $\fa_{k} = \fq_{1, k} \cap \dots \cap \fq_{t, k}$
 of $\fa_{k}\/.$ By \ref{ex.2}(ii), for each $n \in \nn$,
$$
\fa_{k+n} =  \fa_k^{[p^n]}= \fq_{1, k}^{[p^n]} \cap \dots \cap
\fq_{t, k}^{[p^n]}
$$
is a minimal primary decomposition. Since every ideal of $R$ is
$F$-closed, it follows from \ref{pd.0}(iii) and \ref{pd.1a}(iii)
that there is, for each $j = 1, \ldots, t$, an $f$-sequence
$(\fq_{j, n})_{n \in \nn}$ of $\sqrt{\fq_{j, k}}$-primary ideals
of $R$ with
$$
\fq_{j, n}  =  \begin{cases}
             \fq_{j, k}^{[p^{n-k}]} & \text{if $n \geq k$,}   \\
         f^{-(k-n)}(\fq_{j, k}) & \text{if $n < k$.}
\end{cases}
$$
Note that $\fa_{n} = \fq_{1, n} \cap \dots \cap \fq_{t, n}$ for
all $n \in \nn$. It now follows from Theorem \ref{pd.5} that $\fA$
has a primary decomposition.
\end{proof}

We shall extend the result of Proposition \ref{ex.3} to a wider
class of rings.  Recall that a homomorphism of commutative rings
$g : A \lra B$ is said to be {\em pure\/} if, for every $A$-module
$M$, the map $M \otimes_{A}A \lra M \otimes_{A}B$ is injective.
When this is the case, $g$ must be injective, we identify $A$ as a
subring of $B$, and we say that $A$ is a {\em pure subring\/} of
$B$; also, for each ideal $\fa$ of $A$, we have $\fa B \cap A =
\fa$, so that $A$ is Noetherian if $B$ is. Note that if $A$ is a
direct summand of $B$ as an $A$-module, then $A$ is a pure subring
of $B$.

\begin{lem}
\label{ex.3a} Suppose that $R$ is a pure subring of the reduced
commutative Noetherian ring $R'$. Let $\fA$ be a proper ideal of
$R^{\infty}$. Then
\begin{enumerate}
\item $R^{\infty}$ is a pure subring of $R'^{\infty}$;
\item if $\fA R'^{\infty}$ has a primary decomposition as an ideal of $R'^{\infty}$,
then $\fA$ has a primary decomposition in $R^{\infty}$; and
\item if $R'$ is regular, then every finitely generated
proper ideal of $R^{\infty}$ has a primary decomposition.
\end{enumerate}
\end{lem}

\begin{proof} (i) It is straightforward to use \cite[Theorem 7.13]{HM} to show that
$R^{\infty}$ is a pure subring of $R'^{\infty}$.

(ii) This is now immediate, since it follows from part (i) and the
comments preceding the statement of the lemma that $\fA = \fA
R'^{\infty} \cap R^{\infty}$.

(iii) By Proposition \ref{ex.3}, every finitely generated proper
ideal of $R'^{\infty}$ has a primary decomposition. If $\fA$ is a
finitely generated proper ideal of $R^{\infty}$, then $\fA
R'^{\infty}$ is a finitely generated ideal of $R'^{\infty}$, and
so the claim follows from part (ii).
\end{proof}

\begin{thm}
\label{ex.4} Let $R''$ be a pure subring of a regular (commutative
Noetherian) ring of characteristic $p$, and assume that either $R$
is a purely inseparable extension ring of $R''$ or that $R''$ is a
purely inseparable extension ring of $R$. Let $\fa$ be a proper
ideal of $R$. Then
\begin{enumerate}
\item every finitely generated proper ideal of $R^{\infty}$ has a
primary decomposition;
\item the $f$-sequence $\big((\fa^{[p^n]})^F\big)_{n \in \nn}$
has linear growth of primary decompositions;
\item the set\/ $\bigcup_{n \in \nn}
\ass (\fa^{[p^n]})^F$ is finite;
\end{enumerate}
furthermore, if $R$ has a $p^{m_0}$-weak
test element, for some $m_0 \in \nn$, then, in addition,
\begin{enumerate}
\setcounter{enumi}{3}
\item $\fa^*S^{-1}R = (\fa S^{-1}R)^*$ for every multiplicatively closed subset
$S$ of $R$;
\item the ideal of $R^{\infty}$ to which the $f$-sequence
$\big((\fa^{[p^n]})^*\big)_{n \in \nn}$ corresponds has a primary
decomposition, and so the $f$-sequence
$\big((\fa^{[p^n]})^*\big)_{n \in \nn}$ has linear growth of
primary decompositions; and
\item $\bigcup_{n \in \nn} \ass (\fa^{[p^n]})^* \subseteq
\bigcup_{n \in \nn} \ass  (\fa^{[p^n]})^F$, so that the set\/
$\bigcup_{n \in \nn} \ass (\fa^{[p^n]})^*$ is finite.
\end{enumerate}
\end{thm}

\begin{proof} (i) Since $R^{\infty} \cong R''^{\infty}$ (by
\ref{ex.6}) and $R''$ is Noetherian, we can assume that $R$ itself
is a pure subring of a regular (commutative Noetherian) ring of
characteristic $p$, and then the result follows from Lemma
\ref{ex.3a} because a regular ring is reduced.

All the other claims now follow from Theorem \ref{lg.3a}, because
the $f$-sequence $\big((\fa^{[p^n]})^F\big)_{n \in \nn}$
corresponds to the finitely generated proper ideal $\fa
R^{\infty}$ of $R^{\infty}$ (by Theorem \ref{ex.1}), and so has a
primary decomposition by part (i); note that $\fa^{[p^n]}
\subseteq (\fa^{[p^n]})^F \subseteq(\fa^{[p^n]})^*$ for all $n \in
\nn$.
\end{proof}

Note that a ring $R$ satisfying the hypotheses of Theorem \ref{ex.4}
need not have the property that every ideal is $F$-closed.
As an example, consider $K[X^2, X^3]
\subseteq K[X]$ where $K$ is a field of characteristic $2$ and $X$
is an indeterminate. The regular ring
$K[X]$ is a purely inseparable extension of $K[X^2, X^3]$.
The ideal $\fa = (X^2)K[X^2, X^3]$ does not
contain $X^3\/,$ but
$$
 (X^3)^2 = X^6 \in (X^4)K[X^2,X^3] = ((X^2)K[X^2, X^3])^{[2]}\/.
$$
Therefore $\fa$ is not $F$-closed.

\vspace{0.15in}

Our next two results are concerned with the case where $R$ is a
domain and involve the absolute integral closure $R^+$ of $R$; it
should be noted that $R^{\infty}$ can be identified with the
subring $$\{ \sigma \in R^+: \mbox{~there exists~} n \in \nn
\mbox{~such that~} \sigma^{p^n} \in R\}.$$ It is known that, if
$R$ is an excellent regular local ring, then $R^+$ is flat over
$R$: see \cite[Theorem 9.1]{Hunek96} and \cite{HocHun95}.

\begin{prop}
\label{me.2} Suppose that $R$ is a domain and that $R^+$ is
flat over $R$.

\begin{enumerate}
\item Let $\fq$ be a $\fp$-primary ideal of $R$, and let $\fP$ be
the prime ideal of $R^{\infty}$ corresponding to $\fp$. Then $\fq
R^+ \cap R^{\infty}$ is a $\fP$-primary ideal of $R^{\infty}$.
\item Let $\fa$ be a proper ideal of $R$, suppose that
$$
\fa = \fq_1 \cap \ldots \cap \fq_t, \quad \mbox{~where~} \fq_i
\mbox{~is~} \fp_i\mbox{-primary for~} i = 1, \ldots, t,
$$
is a primary decomposition, and let $\fP_1, \ldots, \fP_t$ be the prime ideals
of $R^{\infty}$ corresponding to  $\fp_1, \ldots, \fp_t$ respectively. Then
$$
\fa R^+ \cap R^{\infty} = \bigcap_{i=1}^t (\fq_i R^+ \cap
R^{\infty}), \quad \mbox{~where~} \fq_i R^+ \cap R^{\infty}
\mbox{~is~}\/ \fP_i\mbox{-primary for~} i = 1, \ldots, t,
$$
is a primary decomposition in $R^{\infty}$.
\end{enumerate}
\end{prop}

\begin{proof} (i) Let $\rho \in R^{\infty}$. Then $\rho = r^{1/p^n}$ for some
$r \in R$ and $n \in \nn$.

If $\rho \not\in \fP$, then $r \not\in \fp$
and multiplication by $r$ provides a monomorphism of the $R$-module $R/\fq$.
Since $R^+$ is flat over $R$, it follows that
multiplication by $r$ provides a monomorphism of the $R$-module $(R/\fq) \otimes_R
R^+$; the latter module is isomorphic to $R^+/\fq R^+$, and this has
a submodule isomorphic to $R^{\infty}/ \fq R^+ \cap R^{\infty}$. Hence
multiplication by $r$ provides a monomorphism of the $R^{\infty}$-module
$R^{\infty}/ \fq R^+ \cap R^{\infty}$, and so
multiplication by $r^{1/p^n} = \rho$ also
provides a monomorphism of this $R^{\infty}$-module.

On the other hand, if $\rho \in \fP$, then $r \in \fp$, and some
power of $r$ belongs to $\fq$; hence some power of $r^{1/p^n} =
\rho$ belongs to $\fq R^+ \cap R^{\infty}$. It follows that $\fq
R^+ \cap R^{\infty}$ is a $\fP$-primary ideal of $R^{\infty}$.

(ii) Since $R^+$ is flat over $R$, we have $\fa R^+ =
\bigcap_{i=1}^t \fq_i R^+$, and so $\fa R^+ \cap R^{\infty} =
\bigcap_{i=1}^t (\fq_i R^+ \cap R^{\infty})$. Also, it follows
from part (i) that $\fq_i R^+ \cap R^{\infty}$ is $\fP_i$-primary
for all $i = 1, \ldots, t$.
\end{proof}

\begin{cor}
\label{me.3} Suppose that $R$ is a subring of an excellent regular
local ring $R'$ and that $R'$ is integral over $R$. Let $\fa$ be a
proper ideal of $R$. Then
\begin{enumerate}
\item the ideal of $R^{\infty}$ to which the
$f$-sequence $\big((\fa^{[p^n]})^+\big)_{n \in \nn}$ corresponds
has a primary decomposition;
\item the $f$-sequence
$\big((\fa^{[p^n]})^+\big)_{n \in \nn}$ has linear growth of
primary decompositions;
\item the set\/ $\bigcup_{n \in \nn}
\ass (\fa^{[p^n]})^+$ is finite;
\end{enumerate}
furthermore, if $R$ has a $p^{m_0}$-weak
test element, for some $m_0 \in \nn$, then, in addition,
\begin{enumerate}
\setcounter{enumi}{3}
\item $\fa^*S^{-1}R = (\fa S^{-1}R)^*$ for every multiplicatively closed subset
$S$ of $R$;
\item the ideal of $R^{\infty}$ to which the $f$-sequence
$\big((\fa^{[p^n]})^*\big)_{n \in \nn}$ corresponds has a primary
decomposition, and so the $f$-sequence
$\big((\fa^{[p^n]})^*\big)_{n \in \nn}$ has linear growth of
primary decompositions; and
\item $\bigcup_{n \in \nn} \ass (\fa^{[p^n]})^* \subseteq
\bigcup_{n \in \nn} \ass  (\fa^{[p^n]})^+$, so that the set\/
$\bigcup_{n \in \nn} \ass (\fa^{[p^n]})^*$ is finite.
\end{enumerate}
\end{cor}

\begin{note} Actually, the conclusion of part (iv) can be
established easily, even without the existence of a $p^{m_0}$-weak
test element, by a combination of arguments in the proof of K. E.
Smith \cite[Lemma 2]{Smith01} and in Huneke \cite[p.\
15]{Hunek96}.
\end{note}

\begin{proof} (i)
We have $R'^+ = R^+$; since $R'^+$ is flat over $R'$ (see
\cite[Proposition 9.1]{Hunek96}), it follows from Proposition
\ref{me.2} that the ideal $\fa R^+ \cap R'^{\infty} = \fa R'^+\cap
R'^{\infty} = (\fa R')R'^+ \cap R'^{\infty}$ of $R'^{\infty}$ has
a primary decomposition. Now $R^{\infty}$ is a subring of
$R'^{\infty}$, and so it follows that $\fa R^+ \cap R^{\infty}$,
which is the contraction of $\fa R^+ \cap R'^{\infty}$ to
$R^{\infty}$, has a primary decomposition. However, by Lemma \ref{fs.1}, the
ideal $\fa R^+ \cap R^{\infty}$ of $R^{\infty}$ corresponds to the
$f$-sequence $\big((\fa^{[p^n]})^+\big)_{n \in \nn}$.

All the remaining claims now follow from Theorem \ref{lg.3a},
because $\fa^{[p^n]}
\subseteq (\fa^{[p^n]})^+ \subseteq (\fa^{[p^n]})^*$ for all $n
\in \nn$.
\end{proof}

The next Theorem \ref{me.16} is one of the main results of this section. It should
be noted that many cases of the result of
\ref{me.16}(iv) about localization of tight closure
follow from the result of I. M. Aberbach, Hochster and Huneke in
\cite[Theorem (6.9)]{AHH}; however, the proof below is very different from theirs.
All the other results of \ref{me.16} are believed to be new.

\begin{thm}
\label{me.16} Suppose that $R$ is equidimensional\/
{\rm \cite[p.\ 250]{HM}} and integral over a regular excellent subring
$A$, and let $\fa$ be an ideal of $R$ which is the extension of a proper
ideal $\fb$ of $A$, that is, $\fa = \fb R$. Then

\begin{enumerate}
\item the ideal of $R^{\infty}$ to which the $f$-sequence
$\big((\fa^{[p^n]})^+\big)_{n \in \nn}$ (see\/ {\rm \ref{me.13}} and\/
{\rm \ref{me.14}})
corresponds has a primary decomposition;
\item the $f$-sequence $\big((\fa^{[p^n]})^+\big)_{n \in \nn}$
has linear growth of primary decompositions;
\item the set\/
$\bigcup_{n \in \nn} \ass (\fa^{[p^n]})^+$ is finite;
\end{enumerate}
furthermore, if $R$ has a $p^{m_0}$-weak
test element, for some $m_0 \in \nn$, then, in addition,
\begin{enumerate}
\setcounter{enumi}{3}
\item {\rm (see Aberbach--Hochster--Huneke \cite[Theorem (6.9)]{AHH})}
$\fa^*S^{-1}R = (\fa S^{-1}R)^*$ for every multiplicatively closed
subset $S$ of $R$;
\item the ideal of $R^{\infty}$ to which the $f$-sequence
$\big((\fa^{[p^n]})^*\big)_{n \in \nn}$ corresponds has a primary
decomposition, and so the $f$-sequence
$\big((\fa^{[p^n]})^*\big)_{n \in \nn}$ has linear growth of
primary decompositions; and
\item $\bigcup_{n \in \nn} \ass (\fa^{[p^n]})^* \subseteq
\bigcup_{n \in \nn} \ass  (\fa^{[p^n]})^+$, so that the set\/
$\bigcup_{n \in \nn} \ass (\fa^{[p^n]})^*$ is finite.
\end{enumerate}
\end{thm}

\begin{proof} (i)
Suppose that there are $h$ minimal prime ideals
$\fp_1, \ldots, \fp_h$ of $R$.
Let $i \in \{1, \ldots, h\}$. Since $R$ is equidimensional, we have
$$
\dim A/A \cap \fp_i = \dim R/ \fp_i = \dim R = \dim A
$$
and so $A \cap \fp_i$ is a minimal prime ideal of $A$. Therefore
$A/A \cap \fp_i$ is a direct factor of $A$, and so is itself a
regular excellent domain, and the extension integral domain
$R/\fp_i$ is integral over it. One can now use Lemma \ref{me.14}
to see that it is sufficient to establish part (i) under the
additional hypothesis that $R$ is a domain, and so it is assumed
that that is the case for the remainder of the proof of part (i).
However, in that case, since the class of excellent rings is
closed under localization, it follows from \cite[Theorem
9.1]{Hunek96} that $A^+$ is flat over $A$. Note that $A^+ = R^+$,
and when we identify $R^{\infty}$ and $A^{\infty}$ with subrings
of $A^+ = R^+$ in the natural way, we have $A^{\infty} \subseteq
R^{\infty} \subseteq A^+ = R^+$.

Let $\fq$ be a $\fp$-primary ideal of $A$. We show next that $(\fq R) R^+ \cap R$
has no embedded associated prime ideal. Note that $(\fq R) R^+ = \fq A^+$. Hence
$\left((\fq R) R^+ \cap R\right) \cap A =\fq A^+ \cap A = \fq^+ = \fq$ since $A$ is
regular, so that every ideal of $A$ is tightly closed.
Suppose that $\fr \in \Spec (R)$
is an embedded prime ideal of $(\fq R) R^+ \cap R$, and seek a
contradiction. Then there exists a
minimal prime ideal $\fs \in \ass_R( (\fq R) R^+ \cap R)$ with $\fs \subset
\fr$ (the symbol `$\subset$' is reserved to denote
strict inclusion). Since $R$ is integral over $A$, it follows that $\fp \subseteq
\fs \cap A \subset \fr \cap A$.
This means that there exists $a \in \fr \cap A$ such that
multiplication by $a$ provides a monomorphism of $A/\fq$. Since
$A^+$ is flat over $A$,
multiplication by $a$ provides a monomorphism of $A^+/\fq A^+ = R^+ / (\fq R)R^+$;
hence multiplication by $a$ provides a monomorphism of $R/(\fq R) R^+ \cap R$.
But this is a contradiction, because $\fr$ consists of zerodivisors on
$R/(\fq R) R^+ \cap R$. This shows that
$(\fq R)^+ = (\fq R) R^+ \cap R$ has no embedded prime.

Furthermore, since $A$ is regular, the Frobenius powers of the
$\fp$-primary ideal $\fq$ of $A$ are again $\fp$-primary (see
\ref{ex.2}(ii)), so that, by the above paragraph, $(\fq^{[p^n]} R)
R^+ \cap R = ((\fq R)^{[p^n]})^+$ has no embedded prime (for each
$n \in \nn$). It now follows from \ref{lg.11} that the
$f$-sequence $\big(((\fq R)^{[p^n]})^+\big)_{n \in \nn}$ has
linear growth of primary decompositions, and the ideal of
$R^{\infty}$ to which it corresponds has a primary decomposition.

Now let $\fb = \fq_1 \cap \ldots \cap \fq_t$, where $\fq_i$ is
$\fp_i$-primary for $i = 1, \ldots, t$, be a minimal primary
decomposition of $\fb$ in $A$. Then, by \ref{ex.2}(ii), for each
$n \in \nn$,
$$
\fb^{[p^n]} = \fq_1^{[p^n]} \cap \ldots \cap \fq_t^{[p^n]}, \quad
\mbox{~where~} \fq_i^{[p^n]} \mbox{~is~} \fp_i\mbox{-primary for~}
i = 1, \ldots, t,
$$
is again a minimal primary decomposition. Since $A^+$ is flat over $A$, we have,
again for each $n \in \nn$, that $\fb^{[p^n]} A^+ = \bigcap_{i=1}^t
\fq_i^{[p^n]} A^+$, that is, $(\fb R)^{[p^n]} R^+ = \bigcap_{i=1}^t
(\fq_i R)^{[p^n]} R^+$. Now contract back to $R$, and recall that
$\fa = \fb R$: one obtains that
$$
(\fa^{[p^n]})^+ =  \bigcap_{i=1}^t ((\fq_i R)^{[p^n]})^+ \quad \mbox{~for all~}
n \in \nn.
$$
By the immediately preceding paragraph, the ideal $\fQ_i$ of $R^{\infty}$ to which
the $f$-sequence
$\big(((\fq_i R)^{[p^n]})^+\big)_{n \in \nn}$
corresponds has a primary decomposition (for each $i = 1, \ldots, t$). Therefore
$\bigcap_{i=1}^t \fQ_i$ has a primary decomposition;
however, by \ref{pd.2}(v), the ideal $\bigcap_{i=1}^t \fQ_i$
corresponds to the $f$-sequence
$\big((\fa^{[p^n]})^+\big)_{n \in \nn}$. This completes the proof of part (i).

All the remaining claims now follow from Theorem \ref{lg.3a},
because $\fa^{[p^n]} \subseteq (\fa^{[p^n]})^+ \subseteq
(\fa^{[p^n]})^*$ for all $n \in \nn$, by Lemma \ref{me.14}(i).
\end{proof}

\begin{cor}
\label{me.6} Assume that $R$ is an
equidimensional complete local ring of dimension $d$, and let
$x_1, \ldots, x_d$ form a system of parameters for $R$. Let
$\fa$ denote a proper ideal of $R$ generated by `polynomials' or
`formal power series' in
$x_1, \ldots, x_d$ with coefficients in the prime subfield of $R$ (or even
in a coefficient field $K$ of $R$). Then

\begin{enumerate}
\item the ideal of $R^{\infty}$ to which the $f$-sequence
$\big((\fa^{[p^n]})^+\big)_{n \in \nn}$ corresponds has a primary decomposition;
\item the $f$-sequence $\big((\fa^{[p^n]})^+\big)_{n \in \nn}$
has linear growth of primary decompositions;
\item the set\/
$\bigcup_{n \in \nn} \ass (\fa^{[p^n]})^+$ is finite;
\item $\fa^*S^{-1}R = (\fa S^{-1}R)^*$ for every multiplicatively closed
subset $S$ of $R$;
\item the ideal of $R^{\infty}$ to which the $f$-sequence
$\big((\fa^{[p^n]})^*\big)_{n \in \nn}$ corresponds has a primary
decomposition, and so the $f$-sequence
$\big((\fa^{[p^n]})^*\big)_{n \in \nn}$ has linear growth of
primary decompositions; and
\item $\bigcup_{n \in \nn} \ass (\fa^{[p^n]})^* \subseteq
\bigcup_{n \in \nn} \ass  (\fa^{[p^n]})^+$, so that the set\/
$\bigcup_{n \in \nn} \ass (\fa^{[p^n]})^*$ is finite.
\end{enumerate}
\end{cor}

\begin{proof} By Cohen's structure theorems for complete local rings,
there exists a coefficient field $K$ of $R$, and $R$ is a finite
module over its complete regular local subring $A := K[[x_1,
\ldots, x_d]]$: see \cite[Theorem 29.4]{HM}, for example. Now a
complete local ring is excellent, and so $A$ is an excellent
regular local ring, and $R$ has a $p^{m_0}$-weak test element, for
some $m_0 \in \nn$, by Hochster--Huneke \cite[Theorem
(6.1)(b)]{HocHun94}. All the claims therefore follow from Theorem
\ref{me.16}.
\end{proof}

The final corollary gives a partial extension of the results of
Corollary \ref{me.6}(iv)-(vi) to an equidimensional excellent
local ring.

\begin{cor}
\label{me.17} Assume that $R$ is an
equidimensional excellent local ring of dimension $d$, and let
$x_1, \ldots, x_d$ form a system of parameters for $R$. Let
$\fa$ denote a proper ideal of $R$ generated by `polynomials' in
$x_1, \ldots, x_d$ with coefficients in the prime subfield of $R$. Then

\begin{enumerate}
\item $\fa^*S^{-1}R = (\fa S^{-1}R)^*$ for every multiplicatively closed
subset $S$ of $R$;
\item the ideal of $R^{\infty}$ to which the $f$-sequence
$\big((\fa^{[p^n]})^*\big)_{n \in \nn}$ corresponds has a primary
decomposition, so that the $f$-sequence
$\big((\fa^{[p^n]})^*\big)_{n \in \nn}$ has linear growth of
primary decompositions; and
\item the set\/
$\bigcup_{n \in \nn} \ass (\fa^{[p^n]})^*$ is finite.
\end{enumerate}
\end{cor}

\begin{proof} We first show that we can apply Corollary \ref{me.6} to the
completion $\widehat{R}$ of $R$, by showing that $\widehat{R}$ is also
equidimensional. Let $\widehat{\fp}$ be a minimal prime ideal of $\widehat{R}$,
and let $\fp := \widehat{\fp}\cap R$. Then $\dim R/\fp = d$, and
$\widehat{\fp}/\fp \widehat{R}$ is a minimal prime of the completion
$\widehat{R}/\fp \widehat{R}$ of $R/\fp$. Since $R$ is universally catenary,
it follows from a theorem of L. J. Ratliff (see \cite[Theorem 31.7]{HM}) that
$$
\dim \widehat{R}/\widehat{\fp} = \dim (\widehat{R}/\fp \widehat{R})/
(\widehat{\fp}/\fp \widehat{R}) = \dim \widehat{R}/\fp \widehat{R} = \dim
R/\fp = \dim R = \dim \widehat{R}.
$$
Hence $\widehat{R}$ is equidimensional.

It now follows from Corollary \ref{me.6} that
the ideal of $\widehat{R}^{\infty}$ to which the $f$-sequence
$\big(((\fa\widehat{R})^{[p^n]})^+\big)_{n \in \nn}$ (of ideals of
$\widehat{R}$) corresponds has a primary decomposition. Therefore, by Lemma
\ref{bk.4}, the ideal of $R^{\infty}$ to which the $f$-sequence
$\big(((\fa\widehat{R})^{[p^n]})^+\cap R\big)_{n \in \nn}$ (of ideals of
$R$) corresponds has a primary decomposition; note that $\fa^{[p^n]} \subseteq
((\fa\widehat{R})^{[p^n]})^+\cap R \subseteq ((\fa\widehat{R})^{[p^n]})^*\cap R$
for all $n \in \nn$. But $R$ has
a completely stable
$p^{m_0}$-weak test element, for some $m_0 \in \nn$, by Hochster--Huneke
\cite[Theorem (6.1)(b)]{HocHun94}, and it is easy to use this to see that
$((\fa\widehat{R})^{[p^n]})^*\cap R = (\fa^{[p^n]})^*$ for all $n \in \nn$.
All the claims now follow from Theorem \ref{lg.3a}.
\end{proof}

\bibliographystyle{amsplain}

\end{document}